\documentclass[MSNbibl,number,citesort,seceqn,dvips]{arxbj}
\usepackage{graphicx}

\aid{0}
\volume{19}
\issue{2}
\pubyear{2013}
\firstpage{492}
\lastpage{520}
\doi{10.3150/11-BEJ409}

\makeatletter
\def\bsuffix #1{#1}
\newtheorem{theo}{Theorem}[section]
\newremark{rem}{Remark}[section]
\newtheorem{lemma}{Lemma}[section]
\newremark{anex}{An example}
\newcommand{\cal}{\mathcal}
\newcommand{\eqref}[1]{(\ref{#1})}
\newcommand{\fracf}[2]{({#1})/({#2})}
\newcommand{\fracc}[2]{{#1}/{(#2)}}
\makeatother

\begin{document}
\begin{frontmatter}

\title{Weighted estimation of the dependence function for an extreme-value distribution}
\runtitle{Dependence function for an extreme-value distribution}

\begin{aug}
\author[1]{\fnms{Liang} \snm{Peng}\thanksref{1}\ead[label=e1]{peng@math.gatech.edu}},
\author[2]{\fnms{Linyi} \snm{Qian}\thanksref{2}\ead[label=e2]{lyqian@stat.ecnu.edu.cn}} \and
\author[3]{\fnms{Jingping} \snm{Yang}\corref{}\thanksref{3}\ead[label=e3]{yangjp@math.pku.edu.cn}}
\runauthor{L. Peng, L. Qian and J. Yang} 
\address[1]{School of Mathematics,
Georgia Institute of Technology, Atlanta, GA 30332-0160, USA.\\ \printead{e1}}
\address[2]{School of Finance
and Statistics, East China Normal University, 500 Dongchuan Road,
Shanghai 200241, China. \printead{e2}}
\address[3]{LMEQF and Department of Financial Mathematics,
Center for Statistical Science, Peking University, Beijing, 100871,
China. \printead{e3}}
\end{aug}

\received{\smonth{10} \syear{2010}}
\revised{\smonth{9} \syear{2011}}

%
\begin{abstract}
Bivariate extreme-value distributions have been used in modeling
extremes in environmental sciences and risk management. An important
issue is estimating the dependence function, such as the Pickands
dependence function. Some estimators for the Pickands dependence
function have been studied by assuming that the marginals are known.
Recently, Genest and Segers [\textit{Ann. Statist.} \textbf{37} (2009)
2990--3022]
derived the asymptotic distributions of those proposed estimators
with marginal distributions replaced by the empirical distributions.
In this article, we propose a class of weighted estimators including
those of Genest and Segers (2009) as special cases. We propose a
jackknife empirical likelihood method for constructing
confidence intervals for the Pickands dependence function, which
avoids estimating the complicated asymptotic variance. A~simulation
study demonstrates the effectiveness of our proposed jackknife
empirical likelihood method.
\end{abstract}

%

\end{frontmatter}

\section{Introduction}\label{sec1}
\label{intro} Let $(X_{11},X_{12}),\ldots, (X_{n1},X_{n2})$ be
independent random pairs with common distribution function $F$ and
continuous marginal distributions $F_1(x)=F(x,\infty)$ and
$F_2(y)=F(\infty,y)$. Then the copula of $F$ is defined as
\[
C(x,y)=P\bigl(F_{1}(X_{11})\le x, F_{2}(X_{12})\le y\bigr).
\]
When $C^t(u^{1/t},v^{1/t})=C(u,v)$ holds for all $u,v\in[0, 1]$ and
$t>0$, $C$ is called an extreme value copula and is determined by
the Pickands dependence function, $A$, through the equation
%
\begin{equation}
\label{Pickands}
C(u,v)=\exp\biggl\{\log(uv)A\biggl(\frac{\log(v)}{\log(uv)}\biggr)\biggr\}
\end{equation}
for all
$(u,v)\in(0, 1]^2\setminus\{(1,1)\}$, where $A$ is a convex
function and satisfies $\max(t,1-t)\le A(t)\le1$ for all $0\le t\le
1$ (see Pickands
\cite{r12} and Falk
and Reiss~\cite{r3}).

Write $Y_{ij}=-\log\{F_j(X_{ij})\}$ for $i=1,\ldots,n$, $j=1, 2$ and
\[
H_n(z)=\frac1n\sum_{i=1}^nI\biggl(\frac{Y_{i1}}{Y_{i1}+Y_{i2}}\le z\biggr).\vspace*{-2pt}
\]
We denote $u\wedge v=\min(u, v)$ and $u\vee v=\max(u, v)$
throughout.
Estimators for the
Pickands dependence function $A(t)$ when the marginal distributions
$F_j,j=1,2$ are known have been proposed by Pickands
\cite{r12}, Deheuvels
\cite{r2}, Hall
and Tajvidi~\cite{r7}, and Cap\'era\`a,
Foug\`eres and Genest~\cite{r1}, defined as
\begin{eqnarray*}
A^P(t)&=&\frac n{\sum_{i=1}^n\{Y_{i1}/t\}\wedge\{Y_{i2}/(1-t)\}},
\\[-2pt]
A^D(t)&=&\frac n{\sum_{i=1}^n\{Y_{i1}/t\}\wedge\{
Y_{i2}/(1-t)\}-t\sum_{i=1}^nY_{i1}-(1-t)\sum_{i=1}^nY_{i2}+n},
\\[-2pt]
A^{\mathit{HT}}(t)&=&\frac
n{\sum_{i=1}^n\{\fracf{nY_{i1}}{t\sum_{j=1}^nY_{j1}}\}\wedge\{
\fracc{nY_{i2}}{(1-t)\sum_{j=1}^nY_{j2}}\}},
\\[-2pt]
A^{\mathit{CFG}}(t)&=&\exp\biggl
\{\lambda(t)\int_0^t\frac{H_n(z)-z}{z(1-z)}\,\mathrm{d}z-\bigl(1-\lambda(t)\bigr)
\int_t^1\frac{H_n(z)-z}{z(1-z)}\,\mathrm{d}z\biggr\},\vspace*{-2pt}
\end{eqnarray*}
respectively, where $\lambda(t)\in[0, 1]$ is a weight function and
$A^P(t)$ and $A^D(t)$ are corresponding limits when $t=0$
or $1$. When the marginal distributions are unknown, similar
nonparametric estimators can be obtained by replacing the marginal
distribution $F_j$ by the corresponding empirical distribution
$F_{nj}(x)=\frac1n\sum_{i=1}^nI(X_{ij}\le x)$ or $\hat
F_{nj}(x)=\frac1{n+1}\sum_{i=1}^nI(X_{ij}\le x).$ We denote
these estimators as $\tilde A^P(t), \tilde A^D(t), \tilde A^{\mathit{HT}}(t)$ and $\tilde A^{\mathit{CFG}}(t)$. Recently, Genest
and Segers~\cite{r5} showed
that $\tilde A^P(t), \tilde A^D(t)$ and $\tilde A^{\mathit{HT}}(t)$ have the
same asymptotic distribution as
\[
\hat A^P(t)=\frac n{\sum_{i=1}^n\{Z_{i1}/(1-t)\}\wedge\{ Z_{i2}/t\}}\vspace*{-2pt}
\]
and that $\tilde A^{\mathit{CFG}}(t)$ with $\lambda(t)=t$ has the same asymptotic
distribution as
\[
\hat A^{\mathit{CFG}}(t)=\exp\Biggl\{-\gamma-\frac
1n\sum_{i=1}^n\bigl(Z_{i1}/(1-t)\bigr)\wedge(Z_{i2}/t)\Biggr\},\vspace*{-2pt}
\]
  where
$\gamma=-\int_0^{\infty}\log(x)\mathrm{e}^{-x}\,\mathrm{d}x$ is the Euler constant and
\[
Z_{ij}=-\log\{\hat F_{nj}(X_{ij})\}\qquad\mbox{for }
i=1,\ldots,n, j=1,2.\vspace*{-2pt}
\]
Moreover, Genest
and Segers~\cite{r5} derived
the asymptotic distributions of $\hat A^P(t)$ and $\hat A^{\mathit{CFG}}(t)$
by noting the following important relationship:
\[
\hat A^P(t)=\biggl\{\int_0^1u^{-1}\hat C_n(u^{1-t}, u^t)\,\mathrm{d}u\biggr\}^{-1}\vadjust{\goodbreak}
\]
and
\[
\hat A^{\mathit{CFG}}(t)=\exp\biggl(-\gamma+\int_0^1\{\hat C_n(u^{1-t},
u^t)-I(u>\mathrm{e}^{-1})\}\{u\log(u)\}^{-1}\,\mathrm{d}u \biggr),\vspace*{-2pt}
\]
 where
\[
\hat C_n(u,v)=\frac1n\sum_{i=1}^nI\bigl(\hat F_{n1}(X_{i1})\le u,
\hat F_{n2}(X_{i2})\le v\bigr).\vspace*{-2pt}
\]

In this article, we propose a class of weighted estimators including
$\hat A^P(t)$ and $\hat A^{\mathit{CFG}}(t)$ as special cases. We provide
details in Section~\ref{sec2}. In Section~\ref{sec3} we propose a jackknife empirical
likelihood method to construct confidence intervals for the Pickands
dependence function. Unlike the normal approximation method, this
new method does not need to estimate any additional quantities, such as
asymptotic variance. In Section~\ref{sec4} we report a simulation study
conducted to examine the finite sample behavior of the proposed jackknife
empirical likelihood method. We provide proofs in Section~\ref{sec5}.\vspace*{-2pt}

\section{Weighted estimation}\label{sec2}\vspace*{-2pt}
It follows from (\ref{Pickands}) that
%
\begin{equation}C(u^{1-t},u^t)=u^{A(t)}\qquad\mbox{for all } u\in[0,
1] \mbox{ and all } t\in[0, 1],\label{CA}\vspace*{-2pt}
\end{equation}
which motivates the estimation of $A(t)$ by minimizing the following
weighted distance with respect to $\alpha\ge0$:
\[
\int_0^1\{\hat C_n(u^{1-t}, u^t)-u^{\alpha}\}^2\bar\lambda(u,t)\,\mathrm{d}u,\vspace*{-2pt}
\]
where $\bar\lambda(u,t)\ge0$ is a weight function. Under some
regularity conditions, the foregoing estimator is the solution of
$\alpha$ to the equation
\[
\int_0^1\{\hat C_n(u^{1-t}, u^t)-u^{\alpha}\}u^{\alpha}\{-\log(u)\}\bar
\lambda(u,t)\,\mathrm{d}u=0\vspace*{-2pt}
\]
for $\alpha>0$. This is a special case of the
proposed M-estimators and Z-estimators of B\"ucher,
Dette and Volgushev~\cite{BDV}. Noting that $u^{\alpha}(-\log
u)\bar\lambda(u,t)=C(u^{1-t},u^t)(-\log u)\bar\lambda(u,t)$ and
$\bar\lambda(u,t)$ is any weight function, we propose treating
$C(u^{1-t},u^t)(-\log u)\bar\lambda(u,t)$ as a new weight
function.
This leads us to estimate $A(t)$ by solving the following equation
with respect to $\alpha\ge0$:
%
\begin{equation}\label{est}\int_0^1\{\hat C_n(u^{1-t},
u^t)-u^{\alpha}\}\lambda(u,t)\,\mathrm{d}u=0,\vspace*{-2pt}
\end{equation}
where $\lambda(u,t)\ge0$ is a new weight function. We denote this new
estimator by $\hat A^w_n(t;\lambda)$. When $\lambda(u,t)$ is taken
as $u^{-1}$ or $\{-u\log(u)\}^{-1}$, $\hat A^w_n(t;\lambda)$ becomes
$\hat A^P(t)$ or $\hat A^{\mathit{CFG}}(t)$. Thus, the foregoing class of
estimators includes the known estimators in the literature as
special cases.\vadjust{\goodbreak}

Write $g(\alpha)=\int_0^{1}\{\hat C_n(u^{1-t},
u^t)-u^{\alpha}\}\lambda(u,t)\,\mathrm{d}u.$ Because $u^{\alpha}$ is a
decreasing function of $\alpha$ for each fixed $u\in[0, 1]$,
$g(\alpha)$ is an increasing function of $\alpha$ for each fixed
$t$.
Moreover, $g(0)<0$ and $g(\infty)>0$ when $n$ is sufficiently large. Thus,
(\ref{est}) has a unique solution $\hat A^w_n(t;\lambda)$ for each
large $n$ and $t\in[0, 1]$. Note that this unique solution might not
satisfy that $\max(t, 1-t)\le\hat A^w_n(t;\lambda)\le1$ and $\hat
A^w_n(0;\lambda)=\hat A^w_n(1;\lambda)=1$.

Let $W(u,v)$ denote a tight Gaussian process with mean 0,
covariance\vspace*{-0.8pt}
\[
E\{W(u_1,v_1)W(u_2,v_2)\}=C(u_1\wedge u_2, v_1\wedge
v_2)-C(u_1,v_1)C(u_2,v_2),
\]
 and $W(u,0)=W(0,v)=W(1,1)=0$ for all
$u,v\in[0, 1]$. The asymptotic distribution for the proposed
estimator $\hat A^w_n(t;\lambda)$ is given in the following theorem.

%
\begin{theo}\label{asm}
Suppose that $\frac{\partial^2}{\partial u^2}C(u,v)$,
$\frac{\partial^2}{\partial v^2}C(u,v)$ and
$\frac{\partial^2}{\partial u\,\partial v}C(u,v)$ are defined and
continuous on the sets ${\cal F}_1=\{(u,v)\dvt
0<u<1 \mbox{ and } 0\le v\le1\}$, ${\cal F}_2=\{(u,v)\dvt 0\le
u\le1, 0<v<1\}$ and ${\cal F}_3=\{(u,v)\dvt 0< u< 1, 0<v<1\}$,
respectively.
Also assume that for each fixed $t\in[0, 1]$ the
function $\lambda(u,t)\ge0$ is continuous and not equal to
0 as a function of $u\in(0, 1)$. Furthermore, assume that\vspace*{-0.8pt}
\[
\cases{\displaystyle
\biggl|\frac{\partial^2}{\partial u^2}C(u,v)\biggr|\le\frac{M}{u(1-u)}&\quad
 for   $(u,v)\in{\cal F}_1$,\vspace*{2pt}\cr\displaystyle
\biggl|\frac{\partial^2}{\partial v^2}C(u,v)\biggr|\le\frac{M}{v(1-v)}&\quad
for  $(u,v)\in{\cal F}_2$,\vspace*{2pt}\cr\displaystyle
\biggl|\frac{\partial^2}{\partial u\,\partial v}C(u,v)\biggr|\le\frac
M{u(1-u)}\wedge\frac M{v(1-v)}& \quad for   $(u,v)\in{\cal
F}_3$
}
\]
for some constant $M>0$, $A'(t)$ is continuous on $[0, 1]$, and
there exist $\delta_1>0$ and $\delta_2\in[0, 1/2)$ such that\vspace*{-0.8pt}
%
\begin{equation}
\label{Th1-1}
\cases{\displaystyle
\sup_{0\le t\le1}\sqrt n\int_0^{(n+1)^{-1/((1-t)\vee
t)}}u^{1/2}\lambda(u,t)\,\mathrm{d}u\to0,\vspace*{2pt}\cr\displaystyle
\sup_{0\le t\le1}\sqrt n\int_{(\fracc{n}{n+1})^{1/((1-t)\vee
t)}}^1(1-u)\lambda(u,t)\,\mathrm{d}u\to0,\vspace*{2pt}\cr\displaystyle
\sup_{0\le t\le1}n^{-1/4+\delta_1}\int_{(n+1)^{-1/((1-t)\vee
t)}}^{(\fracc{n}{n+1})^{1/((1-t)\vee t)}} \lambda(u,t)\,\mathrm{d}u\to0,\vspace*{2pt}\cr\displaystyle
\sup_{0\le t\le1}\int_0^1\bigl\{u^{(1-t)\vee t}\bigl(1-u^{(1-t)\vee t}\bigr)\bigr\}
^{\delta_2}\lambda(u,t)\,\mathrm{d}u<\infty,\vspace*{2pt}\cr\displaystyle
\sup_{0\le t\le1}\int_0^1u^{(1-t)\vee t-(1-t)}u^{(1-t)\delta
_2}(1-u^{1-t})^{\delta_2}\lambda(u,t)\,\mathrm{d}u<\infty,\vspace*{2pt}\cr\displaystyle
\sup_{0\le t\le1}\int_0^1u^{(1-t)\vee t-t}u^{t\delta
_2}(1-u^t)^{\delta_2}\lambda(u,t)\,\mathrm{d}u<\infty,\vspace*{2pt}\cr\displaystyle
\sup_{0\le t\le1}\int_{1/2}^1(-\log u)\lambda(u,t)\,\mathrm{d}t<\infty.
}
\end{equation}
%
Then as $n\to\infty$, $\sup_{0\le t\le1}|\hat
A^{w}_n(t;\lambda)-A(t)|=\mathrm{o}_p(1)$. Moreover, suppose that
$\lambda(u,t)$ is continuous in $(0,1)\times[0,1]$ and
\[
|\lambda(u,t_1)-\lambda(u,t_2)|\leq
|t_1-t_2|^{\delta_0}\lambda_0(u), \qquad  t_1,t_2\in[0,1], u\in(0,1)
\]
for some constant $\delta_0>0$ and function
$\lambda_0(u),u\in(0,1)$, where $\lambda_0(u)$ satisfies that
\[
\int^{1/2}_0u^{\alpha}\lambda_0(u)\,\mathrm{d}u<\infty, \qquad
\int^{1}_{1/2}(1-u^{\alpha})\lambda_0(u)\,\mathrm{d}u<\infty
\]
for all $\alpha>0$. Then, as $n\rightarrow\infty$, $\sqrt n\{\hat
A^{w}_n(t;\lambda)-A(t)\}$ converges to $B(t)$ in $C([0,1])$, where
\begin{eqnarray*}
B(t)&=&\biggl\{\int_0^1C(u^{1-t},
u^t)\lambda(u,t)\log(u)\,\mathrm{d}u\biggr\}^{-1}\\
&&{}\times\int_0^1\{W(u^{1-t},
u^t)-C_1(u^{1-t},u^t)W(u^{1-t},1)-C_2(u^{1-t}, u^t)W(1,
u^t)\}\lambda(u,t)\,\mathrm{d}u,
\end{eqnarray*}
$C_1(u,v)=\frac{\partial}{\partial u}C(u,v)$ and
$C_2(u,v)=\frac{\partial}{\partial v}C(u,v)$.
\end{theo}

%
\begin{rem}
Theorem~\ref{asm} still holds when condition (\ref{Th1-1}) is
replaced by
\[
\cases{\displaystyle
\sup_{0\le t\le1}\sqrt n\int_0^{(n+1)^{-2}}u^{1/2}\lambda(u,t)\,\mathrm{d}u\to0,\vspace*{2pt}\cr\displaystyle
\sup_{0\le t\le1}\sqrt n\int_{(\fracc n{n+1})^2}^1 (1-u)\lambda(u,t)\,\mathrm{d}u\to0,\vspace*{2pt}\cr\displaystyle
\sup_{0\le t\le1}n^{-1/4+\delta_1}\int^{(\fracc
{n}{n+1})^2}_{(n+1)^{-2}}\lambda(u,t)\,\mathrm{d}u\to0,\vspace*{2pt}\cr\displaystyle
\sup_{0\le t\le1}\int_0^1u^{\delta_2/2}(1-u)^{\delta_2}\lambda(u,t)\,\mathrm{d}u<\infty
}
\]
for some $\delta_1>0$ and $\delta_2\in[0, 1/2)$. This follows from
the proof of Theorem~\ref{asm} by replacing
$\int_0^{(n+1)^{-1/((1-t)\vee t)}}$,
$\int_{(n+1)^{-1/((1-t)\vee t)}}^{(\fracc n{n+1})^{1/((1-t)\vee t)}}$
and $\int_{(\fracc n{n+1})^{1/((1-t)\vee t)}}^1$ in (\ref{add1})
by $\int_0^{(n+1)^{-2}}$,
$\int_{(n+1)^{-2}}^{(\fracc n{n+1})^2}$
and $\int_{(\fracc n{n+1})^2}^1$, respectively.
\end{rem}

%
\begin{rem}
A common approach to choosing $\lambda(u,t)$ is to minimize the
asymptotic variance
of $\hat A^w_n(t;\lambda)$. This is difficult to do analytically.
Linear combinations of some known
estimators can be considered instead. For example, suppose that the
weight functions $\lambda_1(u),\ldots,\lambda_q(u)$ give the
corresponding estimators $\hat
A^w_{n,1}(t),\ldots, \hat A^w_{n,q}(t)$. Define the class of new
weight functions as
\[
{\cal
F}_0=\Biggl\{\lambda(u,t)\dvt \lambda(u,t)=\sum_{i=1}^qa_i(t)\lambda_i(u),
a_1(t)\ge0,\ldots, a_q(t)\ge0, \sum_{i=1}^qa_i(t)=1\Biggr\}.
\]
Then one
can choose $a_i's$ to minimize the asymptotic variance of $\hat
A^w_n(t;\lambda)$ in this class ${\cal F}_0$, which results in
explicit formulas for $a_i's$.
\end{rem}

\begin{figure}

\includegraphics{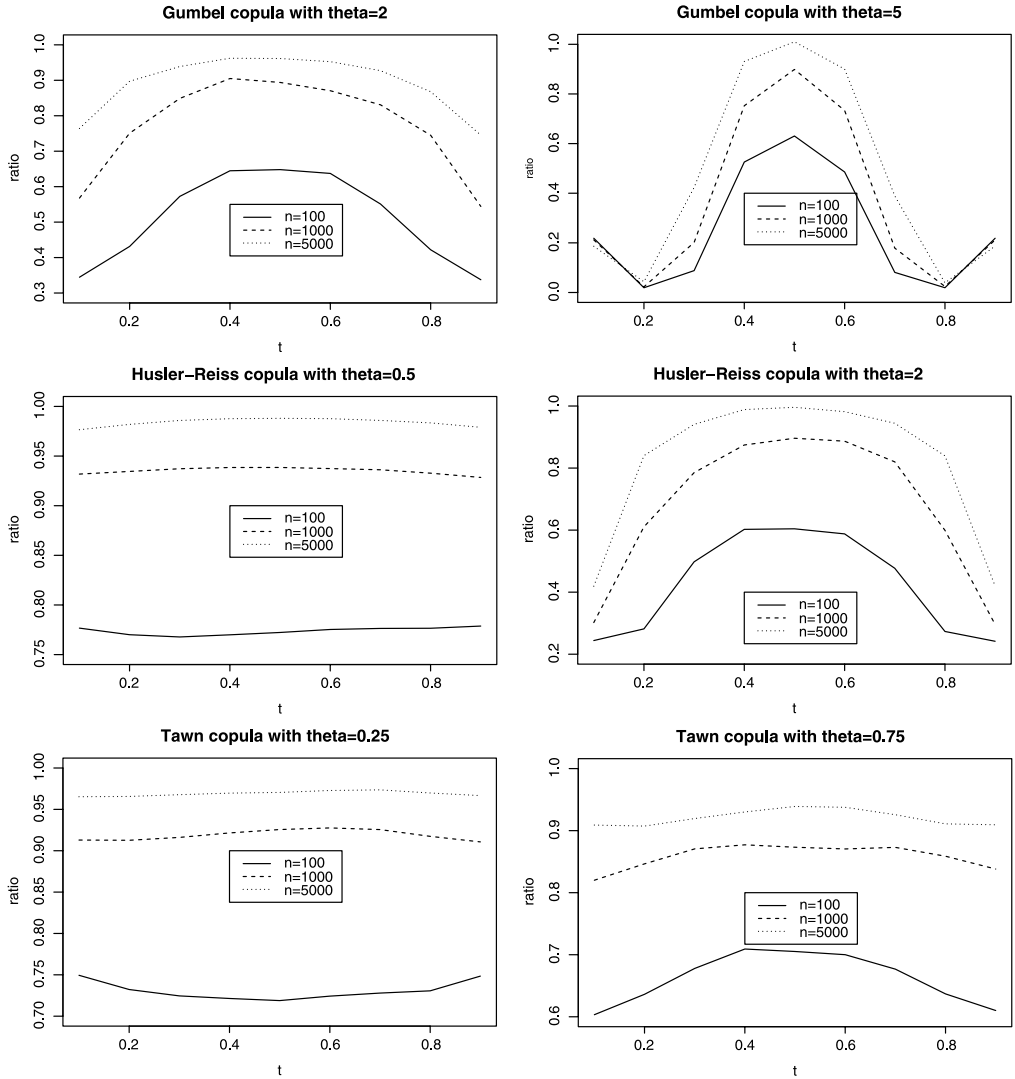}

\caption{Ratios of the mean squared error of the new
estimator $\hat A^w_n(t)$ to that of $\hat A^{\mathit{CFG}}(t)$ for
$t=0.1,0.2,\ldots,0.9$.}\label{ratio}
\end{figure}

\begin{anex*} Assume that $\lambda(u,t)=u^{-1}(-\log
u)^{-q(t)}$ for some $q(t)\in[0, 1]$. Then $\hat A^P(t)$ and $\hat
A^{\mathit{CFG}}(t)$ correspond to $q(t)=0$ and $q(t)=1$, respectively. When
$q(t)<1$, we can write
\begin{eqnarray*}
&&\int_0^1\{\hat C_n(u^{1-t},u^t)-u^{\theta}\}\lambda(u,t)\,\mathrm{d}u\\
&& \quad =-\frac1{1-q(t)}\int_0^1\{\hat C_n(u^{1-t},u^t)-u^{\theta}\}\,\mathrm{d}(-\log u)^{1-q(t)}\\
&& \quad =\frac1{1-q(t)}\int_0^1(-\log u)^{1-q(t)}\,\mathrm{d}\bigl(\hat
C_n(u^{1-t},u^t)-u^{\theta}\bigr)\\
&& \quad =\frac1{1-q(t)}\int_0^1(-\log u)^{1-q(t)}\,\mathrm{d}\hat
C_n(u^{1-t},u^t)-\frac{\theta^{q(t)-1}}{1-q(t)}\int^{\infty
}_0u^{1-q}\mathrm{e}^{-u}\,\mathrm{d}u\\
&& \quad =\frac1{1-q(t)}\Biggl\{\frac1n\sum_{i=1}^n\biggl\{\frac{Z_{i1}}{1-t}\wedge\frac
{Z_{i2}}{t}\biggr\}^{1-q(t)}-\theta^{q(t)-1}\Gamma\bigl(2-q(t)\bigr)\Biggr\},
\end{eqnarray*}
where the $Z_{ij}'s$ are as defined in Section~\ref{sec1}. Thus,
\[
\hat A^w_n(t;\lambda)=\exp\Biggl\{-\Biggl(\log\Biggl(\frac1n\sum_{i=1}^n\biggl(\frac
{Z_{i1}}{1-t}\wedge\frac{Z_{i2}}t\biggr)^{1-q(t)}\Biggr)-\log\Gamma
\bigl(2-q(t)\bigr)\Biggr)\Big/\bigl(1-q(t)\bigr)\Biggr\}
\]
for $0\le q(t)\le1$. Note that when $q(t)=1$, the foregoing expression
is defined as the limit, which becomes the same as $\hat
A^{\mathit{CFG}}(t)$. In particular, we propose to choose $q(t)=\min\{\hat
A^{\mathit{CFG}}(t),1\}$ and denote the resulting estimator by $\hat
A^w_n(t).$ To compare this new estimator with $\hat A^{\mathit{CFG}}(t)$, we
draw $1000$ random samples with size $n=100, 1000, 5000$ from
a Gumbel copula with $A(t)=\{t^{\theta}+(1-t)^{\theta}\}^{1/\theta}$,
a H\"usler--Reiss copula with $A(t)=(1-t)\Phi(\theta+\frac
1{2\theta}\log\frac{1-t}t)+t\Phi(\theta+\frac
1{2\theta}\log\frac{t}{1-t})$, and a Tawn copula with $A(t)=1-\theta
t+\theta t^2$, where $\Phi(x)$ denotes the distribution function of
$N(0,1)$. Figure~\ref{ratio} plots the ratios of the mean
squared error of $\hat A^w_n(t)$ to the mean squared error of $\hat
A^{\mathit{CFG}}(t)$ for
$t=0.1,0.2,\ldots,0.9$, and shows that the new estimator has a
smaller mean squared error than $\hat A^{\mathit{CFG}}(t)$ in all of the cases
considered.
\end{anex*}

\section{Jackknife empirical likelihood method}\label{sec3}
In this section, we consider interval estimation for the Pickands
dependence function $A(t)$, which plays an important role in risk
management since one may be concerned with interval estimation for
$C(u,v)$ at some particular values of $u$ and $v$. Note that an
interval for $A(t)$ can be easily transformed to an interval for a
monotone function of $A(t)$. Moreover, these two intervals have the
same coverage probability, but different interval lengths. Because the
upper tail dependence coefficient can be written as a monotone
function of $A(1/2)$, an interval can be constructed via an interval
for $A(1/2)$.

%

An obvious approach to constructing an interval for $A(t)$ is to use
the normal approximation method based on any one of the estimators
for $A(t)$. Because the asymptotic distribution of any one of the
estimators for $A(t)$ depends on its derivative $A'(t)$, the normal
approximation method requires estimating $A'(t)$ first. In an
alternative approach to constructing confidence intervals, the empirical
likelihood method has been extended and applied in various fields
since Owen
\cite{r9,r10} introduced it for construction of a confidence
interval/region for a mean. (See Owen
\cite{r11} for an overview.)
An important feature of the empirical likelihood method is its property
of self-studentization, which avoids estimating the
asymptotic variance explicitly. A general approach to formulating the
empirical likelihood function is based on estimating equations, as in
Qin
and Lawless~\cite{r13}.

Because our proposed weighted estimator is defined as the solution to
equation (\ref{est}), the method of Qin
and Lawless~\cite{r13} may be applied directly by defining the
empirical likelihood function as
\begin{eqnarray*}
&&\sup\Biggl\{\prod_{i=1}^n(np_i)\dvt p_1\ge0,\ldots,p_n\ge0, \sum
_{i=1}^np_i=1,\\
&&\hphantom{\sup\Biggl\{} \sum_{i=1}^np_i\int_0^{1}\bigl\{I\bigl(\hat F_{n1}(X_{i1})\le u^{1-t}, \hat
F_{n2}(X_{i2})\le u^t\bigr)-u^{\theta}\bigr\}
\lambda(u,t)\,\mathrm{d}u=0 \Biggr\}.
\end{eqnarray*}
However, this method cannot catch the variation
introduced by the marginal empirical distributions. In other words,
the limit is no longer
a chi-squared distribution.
In general, the nonlinear functional must be linearized by introducing
some link variables before the profile empirical likelihood method is used.
(See Chen, Peng and Zhao~\cite{cpz} for details on applying the
profile empirical likelihood method to copulas.) Unfortunately, this
linearization idea is not applicable to the estimation of $A(t)$.
Recently, Jing,
Yuan and Zhou~\cite{r8} proposed a so-called ``jackknife
empirical likelihood'' method to construct confidence intervals for
U-statistics. More specifically, these authors proposed applying the
empirical likelihood method to jackknife samples, which could result in
a chi-squared limit.
Motivated by Gong,
Peng and Qi~\cite{r6}'s study of the use of a smoothed jackknife
empirical likelihood method to construct a confidence interval for a
receiver operating characteristic curve, one needs to work with a
smoothed version of the left-hand side of
(\ref{est}). The reason for smoothing is to separate marginals from the
copula estimator when expanding the jackknife empirical likelihood
ratio. In this work, we used the smoothed empirical copula of
Fermanian,
Radulovi\'c and Wegkamp~\cite{r4}, defined as
\[
\hat C_n^s(u^{1-t},u^t)=\frac1n\sum_{i=1}^nK\biggl(\frac{u-\hat
F_{n1}^{1/(1-t)}(X_{i1})}h\biggr)
K\biggl(\frac{u-\hat F_{n2}^{1/t}(X_{i2})}h\biggr),
\]
where $K(x)=\int_{-\infty}^xk(s)\,\mathrm{d}s$, $k$ is a symmetric density
function with support $[-1, 1]$, and $h=h(n)>0$ is a bandwidth.
Based on this smoothed estimation, a jackknife empirical likelihood
function can be constructed as follows. Put $\hat F_{nj,-i}(x)=\frac
1{n}\sum_{l=1, l\neq i}^nI(X_{lj}\le
x)$ for $j=1,2$ and $i=1,\ldots,n$,
\[
\hat C_{n,-i}^s(u^{1-t},u^t)=\frac1{n-1}\sum_{j=1,j\neq
i}^nK\biggl(\frac{u-\hat F_{n1,-i}^{1/(1-t)}(X_{j1})}h\biggr)K\biggl(\frac{u-\hat
F_{n2,-i}^{1/t}(X_{j2})}{h}\biggr)
\]
for $i=1,\ldots,n$, and define
the jackknife sample as
\[
\hat V_i(u, t)=n\hat C_n^s(u^{1-t}, u^t)-(n-1) \hat C_{n,-i}^s(u^{1-t}, u^t)
\]
for $i=1,\ldots,n$.

We next apply the empirical likelihood method
based on estimating equations of Qin
and Lawless~\cite{r13} to the foregoing jackknife sample. This gives
the jackknife empirical likelihood
function for $\theta=A(t)$ as
\begin{eqnarray*}
L(\theta)&=&\sup\Biggl\{\prod_{i=1}^n(np_i)\dvt p_1\ge0,\ldots, p_n\ge0,
\sum_{i=1}^np_i=1,\\
&&\hphantom{\sup\Biggl\{} \sum_{i=1}^np_i\int_{a_n}^{1-b_n}\{\hat
V_i(u,t)-u^{\theta}\}\lambda(u,t)\,\mathrm{d}u=0\Biggr\},
\end{eqnarray*}
where
$a_n>0$ and $b_n>0$. Note that we use $\int_{a_n}^{1-b_n}$ instead
of $\int_0^1$ in defining the foregoing jackknife empirical likelihood
function, to control the bias term and allow the
possibility of $\lambda(0,t)=\infty$ and $\lambda(1,t)=\infty$.

By the standard Lagrange multiplier technique, we obtain the log
jackknife empirical likelihood ratio as
\[
l(\theta)=-2\log L(\theta)=2\sum_{i=1}^n\log\{1+\beta Q_i(\theta)\},
\]
where
\[
Q_i(\theta)=\int_{a_n}^{1-b_n}\{\hat V_i(u,t)-u^{\theta}\}\lambda(u,t)\,\mathrm{d}u
\]
and $\beta=\beta(\theta)$ satisfies
%
\begin{equation}
\label{Langrage} \frac1n\sum_{i=1}^n\frac{Q_i(\theta)}{1+\beta
Q_i(\theta)}=0.
\end{equation}

%
\begin{theo}\label{jelm}
Suppose that $\frac{\partial^2}{\partial
u^2}C(u,v)$, $\frac{\partial^2}{\partial v^2}C(u,v)$
and $\frac{\partial^2}{\partial u\,\partial v}C(u,v)$ are defined and
continuous on the set ${\cal F}_3=\{(u,v),
0<u<1 \mbox{ and } 0<v<1\}$
and
\begin{eqnarray*}
\biggl|\frac{\partial^2}{\partial u^2}C(u,v)\biggr|&\le&\frac{M}{u(1-u)},\qquad\biggl|\frac
{\partial^2}{\partial v^2}C(u,v)\biggr|\le\frac{M}{v(1-v)},\\
\biggl|\frac
{\partial^2}{\partial u\,\partial v}C(u,v)\biggr|&\le&\frac M{u(1-u)}\wedge\frac M{v(1-v)}
\end{eqnarray*}
for $(u,v)\in{\cal F}_3$ and some constant $M>0$. Let $t$ denote a
fixed point in $(0, 1)$. Assume that the function $\lambda(u,t)\ge
0$ is continuous and not identical to 0 as a function of $u\in
(0, 1)$,
$A'(s)$ is continuous at $s=t$, and
%
\begin{equation}
\label{Th2-1}
\cases{\displaystyle
h=h(n)\to0,\qquad nh\to\infty,\vspace*{2pt}\cr\displaystyle
a_n\to0,\qquad b_n\to0,\qquad h/a_n\to0,\qquad h/b_n\to0,\vspace*{2pt}\cr\displaystyle
n^{-1/4+\delta_1}\int_{a_n}^{1-b_n}\lambda(u,t)\,\mathrm{d}u\to0\qquad\mbox
{for some } \delta_1>0,\vspace*{2pt}\cr\displaystyle
\int_0^1 u^{\delta_2}(1-u^{\delta_2})\lambda(u,t)\,\mathrm{d}u<\infty\qquad
\mbox{for some } \delta_2\in[0, 1/2),\vspace*{2pt}\cr\displaystyle
\sqrt nh^2\int_{a_n}^{1-b_n}u^{-3/2}\lambda(u,t)\,\mathrm{d}u\to0,\vspace*{2pt}\cr\displaystyle
\sqrt nh^2\int_{a_n}^{1-b_n}\{\log u\}^{-1}u^{-3/2}\lambda(u,t)\,\mathrm{d}u\to0,\vspace*{2pt}\cr\displaystyle
\frac1{\sqrt nh}\int_{a_n}^{1-b_n}u^{-1}\lambda(u,t)\,\mathrm{d}u\to0,\vspace*{2pt}\cr\displaystyle
n^{-3/2}\int_{a_n}^{1-b_n}u^{-2}\lambda(u,t)\,\mathrm{d}u\to0}
\end{equation}
as $n\to\infty$. Then
$l(A_0(t))\stackrel{d}{\to}\chi^2(1)$ as $n\to\infty$,
where $A_0(t)$ denotes the true value of $A(t)$.\vspace*{-1.6pt}
\end{theo}

For any fixed $t\in(0, 1)$, based on the foregoing theorem, a jackknife
empirical likelihood confidence interval for $A_0(t)$ with level
$\gamma_0$ can be constructed as
\[
I_{\gamma_0}(t)=\{\theta\dvt l(\theta)\le\chi^2_{\gamma_0}\},
\]
where $\chi^2_{\gamma_0}$ is the $\gamma_0$ quantile of $\chi^2(1)$, as follows:\vspace*{-1.6pt}
%
\begin{rem}
(i) When $\lambda(u,t)=\{-u\log u\}^{-1}$, we have $\sup_{0\leq u\leq
1}\lambda(u,t)=\infty$. We can choose
\[
a_n=d_1n^{-a},\qquad b_n=d_2n^{-b},\qquad h=d_3n^{-1/3}
\]
for some $d_1,d_2,d_3>0$, $0<a<1/9$, and $0<b<1/6$.
{\smallskipamount=0pt
\begin{longlist}[(iii)]
\item[(ii)] When $\sup_{0\le u\le1}\lambda(u,t)<\infty$, we can choose
\[
a_n=d_1n^{-a}, \qquad b_n=d_2n^{-b},\qquad h=d_3n^{-1/3}
\]
for some $d_1,d_2,d_3>0, b>0$, and $0<a<1/3$. Here we fix the rate
for $h$ because the optimal rate for the bandwidth in smoothing
distribution estimation is $n^{-1/3}$.

\item[(iii)] Theorem~\ref{jelm} still holds when $a_n\to a \in(0, 1/2)$ and
$b_n\to b\in(0, 1/2]$ as $n\to\infty$.\vspace*{-1.6pt}
\end{longlist}}
\end{rem}

\section{Simulation study}\label{sec4}
In this section we examine the finite-sample behavior of the
proposed jackknife empirical likelihood method based on
$\lambda(u,t)=u^{-1}(-\log u)^{-\min\{\hat A^{\mathit{CFG}}(t),1\}}$ in terms
of coverage probability and compare it with the method based on
the asymptotic distribution of $\hat A^{\mathit{CFG}}(t)$. For computing the
coverage probability of the proposed jackknife
empirical likelihood method, we choose
$k(x)=\frac{15}{16}(1-x^2)^2I(|x|\le1)$, $h=0.5 n^{-1/3}$,
$a_n=b_n=0.1$, $\lambda(u,t)=u^{-1}(-\log u)^{-\min\{\hat
A^{\mathit{CFG}}(t),1\}}$ and use the R package ``emplik'' (see Zhou~\cite{z}).
To compute the confidence interval based on the asymptotic
distribution of $\hat A^{\mathit{CFG}}(t)$, we use the multiplier method
proposed by Kojadinovic
and Yan~\cite{ky}. More specifically, we use
eq.~(7) of Kojadinovic
and Yan~\cite{ky}, with $N=500$ and
$\{Z_i^{(k)}\dvt i=1,\ldots,n, k=1,\ldots,N\}$ as independent random
variables from $N(0,1)$, to calculate the critical points of the
asymptotic distribution of $\sqrt n\{\hat A^{\mathit{CFG}}(t)-A(t)\}$. We do
not use a larger $N$, because this multiplier method is
computationally intensive. A comparison study on bootstrap
approximations has been reported by B\"ucher
and Dette~\cite{bd}.\looseness=1

%
\begin{table*}
\def\araystretch{1.1}
\tabcolsep=0pt
\caption{Empirical coverage probabilities are reported for the proposed
jackknife empirical likelihood
confidence interval (JELCI) based on $\lambda(u,t)=u^{-1}(-\log
u)^{-\min\{\hat A^{\mathit{CFG}}(t),1\}}$, and the
confidence interval based on the multiplier method for $\hat
A^{\mathit{CFG}}(t)$ (MCI) with nominal levels $0.9$ and $0.95$}\label{tab1}
\begin{tabular*}{\textwidth}{@{\extracolsep{\fill}}lllll@{}}
\hline
& Level 0.9&Level 0.9 & Level 0.95 &Level 0.95\\
$(n,t,\mbox{Copula},\theta)$&JELCI&MCI&JELCI&MCI\\
\hline
$(100,0.1,\mbox{Gumbel},2)$&0.604&0.276&0.639&0.366\\
$(100,0.1,\mbox{H\"usler-Reiss},0.5)$&0.845&0.566&0.899&0.655\\
$(100,0.1,\mbox{Tawn},0.25)$&0.817&0.571&0.872&0.670\\
[5pt]
$(100,0.5,\mbox{Gumbel},2)$&0.871&0.722&0.941&0.784\\
$(100,0.5,\mbox{H\"usler-Reiss},0.5)$&0.888&0.715&0.941&0.802\\
$(100,0.5,\mbox{Tawn},0.25)$&0.886&0.750&0.941&0.825\\
[5pt]
$(100,0.8,\mbox{Gumbel},2)$&0.841&0.531&0.889&0.599\\
$(100,0.8,\mbox{H\"usler-Reiss},0.5)$&0.889&0.646&0.947&0.758\\
$(100,0.8,\mbox{Tawn},0.25)$&0.884&0.677&0.938&0.758\\
[5pt]
$(1000,0.1,\mbox{Gumbel},2)$&0.888&0.655&0.935&0.740\\
$(1000,0.1,\mbox{H\"usler-Reiss},0.5)$&0.892&0.813&0.942&0.883\\
$(1000,0.1,\mbox{Tawn},0.25)$&0.900&0.820&0.957&0.891\\
\hline
\end{tabular*}   \vspace*{6pt}
\end{table*}

We draw $1000$ random samples with size $n=100, 1000$ from the Gumbel
copula, the H\"usler--Reiss copula, and the Tawn copula with Pickands
dependence functions specified at the end of Section~\ref{sec2}. Table
\ref{tab1} reports the coverage probabilities at levels $0.9$ and
$0.95$ for $t=0.1, 0.5, 0.8$. These show that (i) the proposed
jackknife empirical likelihood method gives much more accurate
coverage probabilities than the multiplier method based on the
asymptotic distribution of $\hat A^{\mathit{CFG}}(t)$, and (ii) our proposed
jackknife empirical likelihood method performs poorly for the
boundary case $t=0.1$ when $n=100$, but its performance improves as
$n$ becomes large.\looseness=1

\section{Proofs}\label{sec5}
\vspace*{-\baselineskip}
\begin{pf*}{Proof of Theorem~\ref{asm}}
Define
\[
\alpha_n(u,v)=\sqrt n\Biggl\{\frac1n\sum_{i=1}^nI\bigl(F_1(X_{i1})\le
u, F_2(X_{i2})\le v\bigr)-C(u,v)\Biggr\}.
\]
Then, from Proposition 4.2 of
Segers
\cite{seger} and Theorem G.1 of Genest
and Segers~\cite{r5}, it follows that
\begin{eqnarray*}
&&\sup_{0\le u,v\le1}\bigl|\sqrt n\{\hat C_n(u,v)-C(u,v)\}-\alpha
_n(u,v)+C_1(u,v)\alpha_n(u,1)+C_2(u,v)\alpha_n(1,v)\bigr|\\
&& \quad =\mathrm{O}(n^{-1/4}(\log n)^{1/2}(\log\log n))\qquad
\mbox{a.s.}
\end{eqnarray*}
and
\[
\frac{\alpha_n(u,v)}{(u\wedge v)^{\delta}(1-u\wedge v)^{\delta
}}\stackrel{D}{\to}\frac{W(u,v)}{(u\wedge v)^{\delta}(1-u\wedge
v)^{\delta}}
\]
in the space $l^{\infty}([0, 1]^2)$ of bounded, real-valued
functions on $[0, 1]^2$ for any $\delta\in[0, 1/2)$, where $W$ is
defined before Theorem~\ref{asm}. By the Skorohod construction,
there exists a probability space carrying $\hat C_n^*, \alpha_n^*,
W^*$ such that
%
\begin{equation}
\label{pf1-0} (\hat C_n^*, \alpha_n^*)\stackrel{d}{=}(\hat C_n,
\alpha_n),\qquad W^*\stackrel{d}{=}W,\vspace*{-16pt}
\end{equation}
%
%
\begin{eqnarray}
\label{pf1-1}
&&\sup_{0\le u,v\le1}\bigl|\sqrt n\{\hat C_n^*(u,v)-C(u,v)\}-\alpha
_n^*(u,v)+C_1(u,v)\alpha_n^*(u,1)+C_2(u,v)\alpha_n^*(1,v)\bigr|\nonumber
\\[-8pt]
\\[-8pt]
&& \quad =\mathrm{O}(n^{-1/4}(\log n)^{1/2}(\log\log n))\qquad
\mbox{a.s.}
\nonumber
\end{eqnarray}
and
%
\begin{equation}
\label{pf1-2}
\sup_{0\le u,v\le1} \biggl|\frac{\alpha_n^*(u,v)}{(u\wedge v)^{\delta
}(1-u\wedge v)^{\delta}}-\frac{W^*(u,v)}{(u\wedge v)^{\delta}(1-u\wedge
v)^{\delta}} \biggr|=\mathrm{o}_p(1).
\end{equation}
Let $\hat A^{w*}_n(t;\lambda)$ denote the solution to
\[
\int_0^1\{\hat C_n^*(u^{1-t},u^t)-u^{\alpha}\}\lambda(u,t)\,\mathrm{d}u=0.
\]
Then (\ref{pf1-0}) implies that
%
\begin{equation}
\label{pf1-2a}
\{\hat A^{w*}_n(t;\lambda)\dvt 0\le t\le1\}\stackrel{d}{=}\{\hat
A^w_n(t;\lambda)\dvt 0\le t\le1\}.
\end{equation}

Write
%
\begin{eqnarray}\label{add1}
&&\int_0^{1}\bigl\{\hat C_n^*(u^{1-t},u^{t})-u^{A(t)}\bigr\}\lambda(u,t)\,\mathrm{d}u\nonumber\\
&& \quad =\int_0^{(n+1)^{-1/((1-t)\vee t)}}\bigl\{-u^{A(t)}\bigr\}\lambda(u,t)\,\mathrm{d}u\nonumber\\
&& \qquad {}+\int_{(n+1)^{-1/((1-t)\vee t)}}^{(\fracc{n}{n+1})^{1/((1-t)\vee t)}}\bigl\{
\hat C_n^*(u^{1-t},u^{t})-u^{A(t)}\bigr\}\lambda(u,t)\,\mathrm{d}u\\
&& \qquad {}+\int_{(\fracc{n}{n+1})^{1/((1-t)\vee t)}}^1\bigl\{1-u^{A(t)}\bigr\}\lambda
(u,t)\,\mathrm{d}u\nonumber\\
&& \quad =:I_1(t)+I_2(t)+I_3(t).\nonumber
\end{eqnarray}
Because $1\ge A(t)\ge(1-t)\vee t\ge1/2$, (\ref{Th1-1}) implies that
$I_1(t)$ and $I_3(t)$ are finite and
%
\begin{equation}
\label{pf1-4}
\cases{\displaystyle
\sup_{0\le t\le1}\sqrt n|I_1(t)|\le\sup_{0\le t\le1}\sqrt n\int
_0^{(n+1)^{-1/((1-t)\vee t)}}u^{1/2}\lambda(u,t)\,\mathrm{d}u
=\mathrm{o}(1),\cr\displaystyle
\sup_{0\le t\le1}\sqrt n|I_3(t)|\le\sup_{0\le t\le1}\sqrt n\int
_{(\fracc{n}{n+1})^{1/((1-t)\vee t)}}^1(1-u)
\lambda(u,t)\,\mathrm{d}u=\mathrm{o}(1).
}
\end{equation}
From the condition
\[
\sup_{0\le t\le
1}\int_0^1\bigl\{u^{(1-t)\vee t}\bigl(1-u^{(1-t)\vee
t}\bigr)\bigr\}^{\delta_2}\lambda(u,t)\,\mathrm{d}u<\infty
\]
in (\ref{Th1-1}) and
(\ref{pf1-2}), it follows that
%
\begin{equation}\label{pf1-6(a)}
\sup_{0\le t\le1}\biggl|\int_0^1 \bigl(\alpha_n^*(u^{1-t},u^t)-W^*(u^{1-t},u^t)
\bigr)\lambda(u,t)\,\mathrm{d}u\biggr|
\stackrel{p}{\to}0.
\end{equation}
By (\ref{Pickands}), we have
\begin{eqnarray*}
 0&\le& C_1(u^{1-t}, u^t)=u^{A(t)-(1-t)}\{A(t)-tA'(t)\}\\
  &\le& u^{(1-t)\vee t-(1-t)}\{A(t)-tA'(t)\}
\end{eqnarray*}
and
\begin{eqnarray*}
 0&\le& C_2(u^{1-t}, u^t)=u^{A(t)-t}\{A(t)+(1-t)A'(t)\}\\
  &\le& u^{(1-t)\vee t-t}\{A(t)+(1-t)A'(t)\}.
\end{eqnarray*}
Because $A(t)$ and $A'(t)$ are bounded on [0,1],
from the conditions
\[
\sup_{0\le t\le1}\int_0^1u^{(1-t)\vee t-(1-t)}u^{(1-t)\delta
_2}(1-u^{1-t})^{\delta_2}\lambda(u,t)\,\mathrm{d}u<\infty
\]
and
\[
\sup_{0\le t\le1}\int_0^1u^{(1-t)\vee t-t}u^{t\delta_2}(1-u^t)^{\delta
_2}\lambda(u,t)\,\mathrm{d}u<\infty
\]
in (\ref{Th1-1}), it follows that
%
\begin{equation}\label{pf1-6(b)}
\cases{\displaystyle
\sup_{0\le t\le1}\biggl|\int_0^1\alpha_n^*(u^{1-t},
1)C_1(u^{1-t},u^t)\lambda(u,t)\,\mathrm{d}u\cr\displaystyle
 \hphantom{\sup_{0\le t\le1}\biggl|}-\int_0^1W^*(u^{1-t},1)C_1(u^{1-t},u^t)
\lambda(u,t)\,\mathrm{d}u\biggr|\stackrel{p}{\to}0,\vspace*{2pt}\cr\displaystyle
\sup_{0\le t\le1}\biggl|\int_0^1\alpha_n^*(1,u^t)C_2(u^{1-t},u^t)\lambda
(u,t)\,\mathrm{d}u\cr\displaystyle
\hphantom{\sup_{0\le t\le1}\biggl|} -\int_0^1 W^*(1,u^t)C_2(u^{1-t},u^t)
\lambda(u,t)\,\mathrm{d}u\biggr|\stackrel{p}{\to}0.
}
\end{equation}
By the condition
\[
\sup_{0\le t\le1}n^{-1/4+\delta_1}\int_{(n+1)^{-1/((1-t)\vee
t)}}^{(\fracc n{n+1})^{1/((1-t)\vee t)}}\lambda(u,t)\,\mathrm{d}u\to0
\]
in (\ref{Th1-1}), (\ref{pf1-1}), (\ref{pf1-6(a)}) and
(\ref{pf1-6(b)}), we have
%
\begin{eqnarray}
\label{pf1-7}
&&\sup_{0\le t\le1}\biggl|\sqrt nI_2(t)-\int_0^1\{
W^*(u^{1-t},u^t)-W^*(u^{1-t},1)C_1(u^{1-t},u^t)\nonumber\\
&&\hphantom{\sup_{0\le t\le1}\biggl|\sqrt nI_2(t)-\int_0^1\{}
{}-W^*(1,u^t)C_2(u^{1-t},u^t)\}\lambda(u,t)\,\mathrm{d}u\biggr|
\\
&& \quad =\mathrm{O}_p\biggl(n^{-1/4}(\log n)^{1/2}(\log\log n)^{1/4}
\sup_{0\le t\le1}\int_{(n+1)^{-1/((1-t)\vee t)}}^{(\fracc
n{n+1})^{1/((1-t)\vee t)}}\lambda(u,t)\,\mathrm{d}u \biggr)+\mathrm{o}_p(1) =\mathrm{o}_p(1).
\nonumber
\end{eqnarray}
By (\ref{pf1-4}) and (\ref{pf1-7}), we have
\begin{eqnarray*}
&&\sup_{0\le t\le1}\biggl|\int_0^1\sqrt n\bigl\{\hat C_n^*(u^{1-t}, u^t)-u^{A(t)}\bigr\}
\lambda(u,t)\,\mathrm{d}u\\
&&\hphantom{\sup_{0\le t\le1}\biggl|}{} -\int_0^1\{W^*(u^{1-t},u^t)-W^*(u^{1-t},1)C_1(u^{1-t},u^t)
 -W^*(1,u^t)C_2(u^{1-t},u^t)\}\lambda(u,t)\,\mathrm{d}u\biggr|\\
 && \quad =\mathrm{o}_p(1),
\end{eqnarray*}
which is equivalent to
%
\begin{eqnarray}
\label{pf1-8(a)}
&&\sup_{0\le t\le1}\biggl|\int_0^1\sqrt n\bigl\{u^{\hat
A^{w*}_n(t;\lambda)}-u^{A(t)}\bigr\}\lambda(u,t)\,\mathrm{d}u\nonumber\\
&&\hphantom{\sup_{0\le t\le1}\biggl|}{} -\int_0^1\{W^*(u^{1-t},u^t)-W^*(u^{1-t},1)C_1(u^{1-t},u^t)
\\
 &&\hphantom{\sup_{0\le t\le1}\biggl|{} -\int_0^1\{}
 {}-W^*(1,u^t)C_2(u^{1-t},u^t)\}\lambda(u,t)\,\mathrm{d}u\biggr|   =\mathrm{o}_p(1).\nonumber
\end{eqnarray}
The foregoing equation shows that as $n\rightarrow\infty$,
%
\begin{equation}
\label{dif} \sup_{0\le t\le1}\biggl|\int_0^1\bigl\{u^{\hat
A^{w*}_n(t;\lambda)}-u^{A(t)}\bigr\}\lambda(u,t)\,\mathrm{d}u\biggr|=\mathrm{o}_p(1),
\end{equation}
which implies that
\begin{eqnarray*}
&&P \bigl(\hat A^{w*}_n(t;\lambda)>4/3 \mbox{ for some } t\in[0, 1] \bigr)\\
&& \quad \leq P \biggl(\sup_{0\le t\le1}\biggl|\int_0^1\bigl\{u^{\hat
A^{w*}_n(t;\lambda)}-u^{A(t)}\bigr\}\lambda(u,t)\,\mathrm{d}u\biggr| \geq\inf_{0\le t\le
1}\int_0^1\bigl(u^{A(t)}-u^{4/3}\bigr)\lambda(u,t)\,\mathrm{d}u \biggr)  \rightarrow0
\end{eqnarray*}
since $1/2\le A(t)\le1$ for all $0\le t\le1$. Similarly,
\begin{eqnarray*}
&&P \bigl(\hat A^{w*}_n(t;\lambda)<1/3 \mbox{ for some } t\in[0, 1] \bigr)\\
&& \quad \leq P \biggl(\sup_{0\le t\le1}\biggl|\int_0^1\bigl\{u^{\hat
A^{w*}_n(t;\lambda)}-u^{A(t)}\bigr\}\lambda(u,t)\,\mathrm{d}u\biggr|   \geq\inf_{0\le t\le
1}\int_0^1\bigl(u^{1/3}-u^{A(t)}\bigr)\lambda(u,t)\,\mathrm{d}u \biggr)\\
&&\hspace*{6.7pt}\rightarrow0\qquad\mbox{as } n\to\infty.
\end{eqnarray*}
Thus,
%
\begin{equation}
\label{pf1-9}
P \bigl(1/3\le\hat A^{w*}_n(t;\lambda)\le4/3 \mbox{ for all } t\in[0,1]
\bigr)\to1.
\end{equation}
By the mean value theorem,
%
\begin{eqnarray}
\label{pf1-10}
&&\int_0^1\bigl\{u^{\hat
A^{w*}_n(t;\lambda)}-u^{A(t)}\bigr\}\lambda(u,t)\,\mathrm{d}u\nonumber\\
&& \quad =\int_0^1 u^{a(u,t)A(t)+(1-a(u,t))\hat A^{w*}_n(t;\lambda)}(\log
u)\lambda(u,t)\,\mathrm{d}u\\
&& \qquad {} \times\bigl(A(t)-\hat A^{w*}_n(t;\lambda)\bigr)\nonumber
\end{eqnarray}
for some $a(u,t)\in[0, 1]$.
Because $1/2\le A(t)\le1$, we have $0<a(u,t)A(t)+(1-a(u,t))\hat
A^{w*}_n(t; \lambda)\le7/3$ when $0<\hat A^{w*}_n(t;\lambda)\le4/3$. Thus,
from (\ref{pf1-9}), it follows that
\begin{eqnarray*}
&&P \biggl(\inf_{0\le t\le1}\int_0^1 u^{a(u,t)A(t)+(1-a(u,t))\hat
A^{w*}_n(t;\lambda)}(-\log u)
\lambda(u,t)\,\mathrm{d}u \\
&&\hphantom{P \biggl(} \ge\sup_{0\le t\le1} \int_0^1u^{7/3}(-\log
u)\lambda(u,t)\,\mathrm{d}u \biggr)\rightarrow1
\end{eqnarray*}
as $n\to\infty$,
which, combined with (\ref{dif}), (\ref{pf1-10}) and (\ref{pf1-0}),
implies that
%
\begin{equation}
\label{pf1-10a}
\sup_{0\le t\le1}|\hat A^{w*}_n(t;\lambda)-A(t)|=\mathrm{o}_p(1).
\end{equation}
Then $\sup_{0\le t\le1}|\hat A^{w}_n(t;\lambda)-A(t)|=\mathrm{o}_p(1)$
follows from \eqref{pf1-10a} and \eqref{pf1-2a}.

We next prove that $\hat{A}_n^w(t;\lambda)$ is continuous for $t\in
[0,1]$. For $t_m, t\in[0,1]$ and $t_m\rightarrow t\in[0,1]$ as
$m\rightarrow\infty$, we have
\begin{eqnarray*}
&&\int^{1/2}_0u^{\hat{A}_n^w(t_m;\lambda)}\lambda(u,t_m)\,\mathrm{d}u+\int
^1_{1/2}\bigl(u^{\hat{A}_n^w(t_m;\lambda)}-1\bigr)\lambda(u,t_m)\,\mathrm{d}u\\
&& \quad =\int^{1/2}_0 \hat C_n(u^{1-t_m},
u^{t_m})\lambda(u,t_m)\,\mathrm{d}u+\int^1_{1/2}\bigl( \hat C_n(u^{1-t_m},
u^{t_m})-1\bigr)\lambda(u,t_m)\,\mathrm{d}u.
\end{eqnarray*}
Note that the function
\[
\int^{1/2}_0 \hat C_n(u^{1-t},
u^t)\lambda(u,t)\,\mathrm{d}u+\int^1_{1/2}\bigl( \hat C_n(u^{1-t},
u^t)-1\bigr)\lambda(u,t)\,\mathrm{d}u
\]
is continuous in $t\in[0,1]$; thus, we have
\begin{eqnarray*}
&&\lim_{m\rightarrow\infty}
\biggl(\int^{1/2}_0u^{\hat{A}_n^w(t_m;\lambda)}\lambda(u,t_m)\,\mathrm{d}u+\int
^1_{1/2}\bigl(u^{\hat{A}_n^w(t_m;\lambda)}-1\bigr)
\lambda(u,t_m)\,\mathrm{d}u \biggr)\\
&& \quad =\int^{1/2}_0 \hat C_n(u^{1-t}, u^t)\lambda(u,t)\,\mathrm{d}u+\int^1_{1/2}\bigl(
\hat C_n(u^{1-t}, u^t)-1\bigr)\lambda(u,t)\,\mathrm{d}u.
\end{eqnarray*}
Because
\[
\int^{1/2}_0u^{\alpha}\lambda(u,t)\,\mathrm{d}u+\int^1_{1/2}(u^{\alpha}-1)\lambda(u,t)\,\mathrm{d}u
\]
is continuous in $t\in[0,1]$ and is monotone in $\alpha$ for each
$t\in[0,1]$, we conclude that
$\hat{A}_n^w(t_m;\lambda)\rightarrow\hat{A}_n^w(t;\lambda)$ as
$m\rightarrow\infty$. Thus, $\hat{A}_n^w(t;\lambda)$ is continuous
in $[0,1]$.

Note that
\[
\sup_{0\le t\le1}\int_0^1\bigl\{u^{(1-t)\vee t}\bigl(1-u^{(1-t)\vee t}\bigr)\bigr\}^{\delta
_2}\lambda(u,t)\,\mathrm{d}u<\infty
\]
for some $\delta_2\in[0, 1/2)$ in (\ref{Th1-1})
implies that $ \int_0^{1/2}u^{\delta_2}\lambda(u,t)\,\mathrm{d}u<\infty.$ Thus, using
\begin{eqnarray*}
&&u^{a(u,t)A(t)+(1-a(u,t))\hat A^{w*}_n(u;\lambda)}(-\log
u)\lambda(u,t)\\
&& \quad =u^{A(t)}(-\log u)\lambda(u,t)u^{(1-a(u,t))(
\hat A^{w*}_n(u;\lambda)-A(t))}\\
&& \quad \leq u^{A(t)}(-\log u)\lambda(u,t)u^{-(1-a(u,t))\sup_{0\le t\le
1}|\hat A^{w*}_n(t;\lambda)-A(t)|},
\end{eqnarray*}
$A(t)\ge1/2$ for all $t\in[0,1]$, (\ref{pf1-10a})
and
\[
0\le u^{-s_1}-1\le\frac{s_1}{s_2}u^{-s_2}\qquad\mbox{for all } u\in
[0, 1] \mbox{ and any fixed } 0<s_1<s_2<1,
\]
we get that
%
\begin{eqnarray}
\label{pf1-12}
&&\sup_{0\le t\le1}\biggl|\int_0^1u^{a(u,t)A(t)+(1-a(u,t))
\hat A^{w*}_n(t;\lambda)}(\log u)\lambda(u,t)\,\mathrm{d}u -\int_0^1u^{A(t)}(\log u)\lambda(u,t)\,\mathrm{d}u\biggr|\nonumber\\
&& \quad \leq\sup_{0\le t\le1}\biggl|\int_0^1u^{A(t)}(-\log
u)\lambda(u,t)\bigl(u^{-(1-a(u,t))\sup_{0\le s\le1}|\hat
A^{w*}_n(s;\lambda)-A(s)|}-1\bigr)\,\mathrm{d}u\biggr|\nonumber\\
&& \quad \le\sup_{0\le t\le1}\biggl (\frac{(1-a(u,t))\sup_{0\le s\le1}|\hat
A^{w*}_n(s;\lambda)-A(s)|}{(1-a)\sup_{0\le s\le1}|\hat
A^{w*}_n(s;\lambda)-A(s)|+(A(t)-\delta_2)/2} \nonumber
\\[-12pt]
\\[-4pt]
&&\hphantom{\quad \le\sup_{0\le t\le1}\biggl (} {}\times\int_0^1u^{-(1-a(u,t))\sup_{0\le s\le1}|\hat A^{w*}_n(s;\lambda
)-A(s)|+(A(t)+\delta_2)/2}(-\log u)\lambda(u,t)\,\mathrm{d}u \biggr)\nonumber\\
&& \quad =\mathrm{o}_p(1)\mathrm{O}_p\biggl(\sup_{0\le t\le1}\int_0^1u^{-(1-a(u,t))\sup_{0\le s\le
1}|\hat A^{w*}_n(s;\lambda)-A(s)|+(A(t)+\delta_2)/2}(-\log u)\lambda
(u,t)\,\mathrm{d}u \biggr)\nonumber\ \ \\
&& \quad =\mathrm{o}_p(1)\mathrm{O}_p\biggl(\sup_{0\le t\le
1}\int_0^{1}u^{\delta_2}(1-u^{\delta_2})\lambda(u,t)\,\mathrm{d}u \biggr)=\mathrm{o}_p(1).\nonumber
\end{eqnarray}
Note that the two processes $\hat A^{w}_n(t;\lambda), B(t)$ are
continuous for $t\in[0,1]$. Thus, from (\ref{pf1-0}),
(\ref{pf1-8(a)}), (\ref{pf1-10}), (\ref{pf1-12}) and \eqref{pf1-2a},
we conclude that $\sqrt n\{\hat A^{w}_n(t;\lambda)-A(t)\}$ converges
to $B(t)$ in $C([0,1])$.
\end{pf*}

Before proving Theorem~\ref{jelm}, we present some lemmas. Throughout, we
assume that $t$ is a given point in $(0, 1)$ and use $\theta_0$
to denote $A_0(t)$.

%
\begin{lemma}\label{lem5.1}
Under the conditions of Theorem~\ref{jelm}, as
$n\rightarrow\infty$, we have
\begin{eqnarray*}
\frac1{\sqrt n}\sum_{i=1}^nQ_i(\theta_0)&\stackrel{d}{\to}&
\int_0^1 \{W(u^{1-t},u^t)-W(u^{1-t},1)C_1(u^{1-t},u^t)\\
&&\hphantom{\int_0^1 \{}{} -W(1,u^t)C_2(u^{1-t},u^t) \}\lambda(u,t)\,\mathrm{d}u.
\end{eqnarray*}
\end{lemma}
\begin{pf}Write
%
\begin{eqnarray}
\label{EXP}
\hat V_i(u,t)&=&K\biggl(\frac{u-\hat F_{n1,-i}^{1/(1-t)}(X_{i1})}h\biggr)K\biggl(\frac
{u-\hat F_{n2,-i}^{1/t}(X_{i2})}h\biggr)\nonumber\\
&&{} +\sum_{j=1}^n\biggl\{K\biggl(\frac{u-\hat F_{n1}^{1/(1-t)}(X_{j1})}h\biggr)K\biggl(\frac
{u-\hat
F_{n2}^{1/t}(X_{j2})}h\biggr)\nonumber
\\[-8pt]
\\[-8pt]
&&\hphantom{+\sum_{j=1}^n\biggl\{}{}-  K\biggl(\frac{u-F_{n1,-i}^{1/(1-t)}(X_{j1})}h\biggr)K\biggl(\frac
{u-F_{n2,-i}^{1/t}(X_{j2})}h\biggr)\biggr\}\nonumber\\
&=:&\hat V_{i1}(u,t)+\hat V_{i2}(u,t)\nonumber
\end{eqnarray}
and
%
\begin{eqnarray}
\label{I1-4}
&&\frac1n\sum_{i=1}^nQ_i(\theta_0)\nonumber\\
&& \quad =n^{-1}\int_{a_n}^{1-b_n}\sum_{i=1}^n\{\hat V_{i1}(u,t)-u^{\theta}\}
\lambda(u,t)\,\mathrm{d}u\nonumber\\
&& \qquad {}+n^{-1}\int_{a_n}^{1-b_n}\sum_{i=1}^n\hat V_{i2}(u,t)\lambda(u,t)\,\mathrm{d}u\nonumber\\
&& \quad =n^{-1}\int_{a_n}^{1-b_n} \Biggl\{\sum_{i=1}^nK\biggl(\frac{u-\hat
F_{n1,-i}^{1/(1-t)}(X_{i1})}h\biggr)
K\biggl(\frac{u-\hat F_{n2,-i}^{1/t}(X_{i2})}h\biggr)-u^{\theta_0} \Biggr\}\lambda(u,t)\,\mathrm{d}u\nonumber\\
&& \qquad {}+n^{-1}\int_{a_n}^{1-b_n}\sum_{i=1}^n\sum_{j=1}^n \biggl\{K\biggl(\frac{u-\hat
F_{n1}^{1/(1-t)}(X_{j1})}h\biggr)
-K\biggl(\frac{u-\hat F_{n1,-i}^{1/(1-t)}(X_{j1})}h\biggr) \biggr\}\nonumber\\
&& \hphantom{{}+n^{-1}\int_{a_n}^{1-b_n}\sum_{i=1}^n\sum_{j=1}^n}\qquad {}\times K\biggl(\frac{u-\hat F_{n2}^{1/t}(X_{j2})}h\biggr)\lambda(u,t)\,\mathrm{d}u\\
&& \qquad {}+n^{-1}\int_{a_n}^{1-b_n}\sum_{i=1}^n\sum_{j=1}^nK\biggl(\frac{u-\hat
F_{n1}^{1/(1-t)}(X_{j1})}h\biggr)\nonumber\\
&& \hphantom{{}+n^{-1}\int_{a_n}^{1-b_n}\sum_{i=1}^n\sum_{j=1}^n}\qquad {} \times\biggl\{K\biggl(\frac{u-\hat F_{n2}^{1/t}(X_{j2})}h\biggr)-K\biggl(\frac{u-\hat
F_{n2,-i}^{1/t}(X_{j2})}h\biggr) \biggr\}\lambda(u,t)\,\mathrm{d}u\nonumber\\
&& \qquad {}+n^{-1}\int_{a_n}^{1-b_n}\sum_{i=1}^n\sum_{j=1}^n \biggl\{K\biggl(\frac{u-\hat
F_{n1,-i}^{1/(1-t)}(X_{j1})}h\biggr)
-K\biggl(\frac{u-\hat F_{n1}^{1/(1-t)}(X_{j1})}h\biggr) \biggr\}\nonumber\\
&& \hphantom{{}+n^{-1}\int_{a_n}^{1-b_n}\sum_{i=1}^n\sum_{j=1}^n}\qquad {} \times\biggl\{K\biggl(\frac{u-F_{n2}^{1/t}(X_{j2})}h\biggr)-K\biggl(\frac{u-\hat
F_{n2,-i}^{1/t}(X_{j2})}h\biggr) \biggr\}\lambda(u,t)\,\mathrm{d}u\nonumber\\
&& \quad =:I_1+\cdots+I_4. \nonumber
\end{eqnarray}
Furthermore, the first term $I_1$ can be expressed as
%
\begin{eqnarray}
\label{I1}
I_1 &=&\int_{a_n}^{1-b_n}\lambda(u,t) \biggl\{\int_0^1\int_0^1\frac
1n\sum_{i=1}^nI\bigl(\hat F_{n1,-i}(X_{i1})\le s_1, \hat
F_{n2,-i}(X_{i2})\le s_2\bigr)h^{-2}\nonumber\\[2pt]
&&\hphantom{\int_{a_n}^{1-b_n}\lambda(u,t) \biggl\{}{}\times
k\biggl(\frac{u-s_1^{1/(1-t)}}h\biggr)k\biggl(\frac
{u-s_2^{1/t}}h\biggr)\,\mathrm{d}s_1^{1/(1-t)}\,\mathrm{d}s_2^{1/t}-u^{\theta_0} \biggr\}\,\mathrm{d}u\\[2pt]
 &=&\int_{a_n}^{1-b_n}\lambda(u,t)\Biggl \{\int_0^1\int_0^1\frac1n\sum
_{i=1}^nI \biggl(\hat F_{n1}(X_{i1})\le\frac n{n+1}\biggl(s_1+\frac1n\biggr), \nonumber\\[2pt]
&& \hphantom{ \int_{a_n}^{1-b_n}\lambda(u,t)\Biggl \{\int_0^1\int_0^1\frac1n\sum
_{i=1}^nI \biggl(}\hat F_{n2}(X_{i2})\le\frac n{n+1}\biggl(s_2+\frac1n\biggr) \biggr)\nonumber\\[2pt]
&&\hphantom{ \int_{a_n}^{1-b_n}\lambda(u,t)\Biggl \{\int_0^1\int_0^1}{}\times
h^{-2}k\biggl(\frac{u-s_1^{1/(1-t)}}h\biggr)\nonumber\\[2pt]
&&\hphantom{ \int_{a_n}^{1-b_n}\lambda(u,t)\Biggl \{\int_0^1\int_0^1}{}\times
k\biggl(\frac
{u-s_2^{1/t}}h\biggr)\,\mathrm{d}s_1^{1/(1-t)}\,\mathrm{d}s_2^{1/t}-u^{\theta_0}
\Biggr\}\,\mathrm{d}u\nonumber\\[2pt]
  &=&
\int_{a_n }^{1-b_n}\int_{-1}^1\int_{-1}^1\lambda(u,t) \biggl\{\hat
C_n\biggl(\frac n{n+1}(u-s_1h)^{1-t}+
\frac1{n+1}, \frac n{n+1}(u-s_2h)^{t}+\frac1{n+1}\biggr)\nonumber\\[2pt]
&&\hphantom{ \int_{a_n }^{1-b_n}\int_{-1}^1\int_{-1}^1\lambda(u,t) \biggl\{}{}
-C\biggl(\frac{n}{n+1}(u-s_1h)^{1-t}+\frac1{n+1},\nonumber\\[2pt]
&&\hphantom{ \int_{a_n }^{1-b_n}\int_{-1}^1\int_{-1}^1\lambda(u,t) \biggl\{{}
-C\biggl(}
\frac n{n+1}(u-s_2h)^t
+\frac1{n+1}\biggr) \biggr\}k(s_1)k(s_2)\,\mathrm{d}s_1\,\mathrm{d}s_2\,\mathrm{d}u\nonumber\\[2pt]
&&   {}+\int_{a_n}^{1-b_n}\int_{-1}^1\int_{-1}^1\lambda(u,t) \biggl\{C\biggl(\frac
{n}{n+1}(u-s_1h)^{1-t}+\frac1{n+1}, \frac{n}{n+1}(u-s_2h)^t+\frac
1{n+1}\biggr)\nonumber\\[2pt]
&&\hphantom{{}+\int_{a_n}^{1-b_n}\int_{-1}^1\int_{-1}^1\lambda(u,t) \biggl\{}   {}-
C(u^{1-t}, u^t) \biggr\}k(s_1)k(s_2)\,\mathrm{d}s_1\,\mathrm{d}s_2\,\mathrm{d}u\nonumber\\[2pt]
  &=:&\mathit{II}_1+\mathit{II}_2.\nonumber
\end{eqnarray}
Because $\sup_{a_n\le u\le1-b_n}(h/u)\le h/a_n\to0$ and
\[
\inf_{a_n\le u\le1-b_n}\min\{(n+1)u^t,(n+1)u^{1-t}\}\ge
(n+1)a_n\to\infty
\]
as $n\to\infty$, we have
%
\begin{eqnarray}
\label{pfL1-1}
&&\sup_{a_n\le u\le1-b_n, -1\le s\le1}\biggl|\frac{\log(u-sh)}{\log
u}-1\biggr|\nonumber
\\[-4pt]
\\[-4pt]
&& \quad \le\sup_{a_n\le u\le1-b_n}\frac{2h/u}{-\log u}\le\frac{2h}{-a_n\log
(a_n)}+\frac{2h}{-(1-b_n)\log(1-b_n)}\to0,
\nonumber
\end{eqnarray}\vspace*{-10pt}
%
%
\begin{equation}
\label{pfL1-1a}
\sup_{a_n\le u\le1-b_n, -1\le s\le1}\biggl|u^{t-1}\biggl\{\frac
n{n+1}(u-sh)^{1-t}+\frac1{n+1}\biggr\}-1\biggr|\to0
\end{equation}
and
%
\begin{equation}
\label{pfL1-1b}
\sup_{a_n\le u\le1-b_n, -1\le s\le1}\biggl|u^{-t}\biggl\{\frac
n{n+1}(u-sh)^{t}+\frac1{n+1}\biggr\}-1\biggr|\to0,
\end{equation}
which, together with (\ref{Pickands}), imply that
%
\begin{equation}
\label{pfL1-2}\hspace*{-10pt}
\cases{\displaystyle
\sup_{a_n\le u\le1-b_n, -1\le s_1,s_2\le1}\biggl|\biggl(\log\biggl\{\frac
n{n+1}(u-s_1h)^{1-t}+\frac1{n+1}\biggr\}\cr\displaystyle
\hphantom{\sup_{a_n\le u\le1-b_n, -1\le s_1,s_2\le1}\biggl|\biggl(}{}+\log\biggl\{\frac n{n+1}(u-s_2h)^t+\frac
1{n+1}\biggr\}\biggr)\big/\log u
 -1\biggr|\to0,\cr\displaystyle
\sup_{a_n\le u\le1-b_n, -1\le s_1,s_2\le1}\biggl|A\biggl( \log\biggl\{\frac
n{n+1}(u-s_2h)^t+\frac1{n+1}\biggr\}\cr\displaystyle\hphantom{\sup_{a_n\le u\le1-b_n, -1\le s_1,s_2\le1}\biggl|A\biggl(}
{}\Big/\biggl(\log\biggl\{\frac
n{n+1}(u-s_1h)^{1-t}+\frac1{n+1}\biggr\}
\cr\displaystyle\hphantom{\sup_{a_n\le u\le1-b_n, -1\le s_1,s_2\le1}\biggl|A\biggl({}\Big/\biggl(}
{}+\log\biggl\{\frac n{n+1}(u-s_2h)^t+\frac
1{n+1}\biggr\}\biggr)\biggr)
 -A(t)\biggr|\to0,\vspace*{2pt}\cr\displaystyle
\sup_{a_n\le u\le1-b_n, -1\le s_1,s_2\le1}\biggl|C\biggl(\frac
{n}{n+1}(u-s_1h)^{1-t}+\frac1{n+1},\cr\displaystyle
\hphantom{\sup_{a_n\le u\le1-b_n, -1\le s_1,s_2\le1}\biggl|C\biggl(}\frac{n}{n+1}(u-s_2h)^t+\frac
1{n+1}\biggr)
 -C(u^{1-t},u^t)\biggr|\to0,\cr\displaystyle
\sup_{a_n\le u\le1-b_n, -1\le s_1,s_2\le1}\biggl|C_1\biggl(\frac
n{n+1}(u-s_1h)^{1-t}+\frac1{n+1},\cr\displaystyle \hphantom{\sup_{a_n\le u\le1-b_n, -1\le s_1,s_2\le1}\biggl|C_1\biggl(}\frac{n}{n+1}(u-s_2h)^t+\frac
1{n+1}\biggr)
 -C_1(u^{1-t},u^t)\biggr|\to0,\cr\displaystyle
\sup_{a_n\le u\le1-b_n, -1\le s_1,s_2\le1}\biggl|C_2\biggl(\frac
{n}{n+1}(u-s_1h)^{1-t}+\frac1{n+1},\cr\displaystyle \hphantom{\sup_{a_n\le u\le1-b_n, -1\le s_1,s_2\le1}\biggl|C_2\biggl(} \frac{n}{n+1}(u-s_2h)^t+\frac
1{n+1}\biggr)
 -C_2(u^{1-t},u^t)\biggr|\to0,\vspace*{1pt}\cr\displaystyle
\sup_{a_n\le u\le1-b_n, -1\le s_1,s_2\le1}\biggl|C_{11}\biggl(\frac
n{n+1}(u-s_1h)^{1-t}+\frac1{n+1},\cr\displaystyle \hphantom{\sup_{a_n\le u\le1-b_n, -1\le s_1,s_2\le1}\biggl|C_{22}\biggl(}\frac{n}{n+1}(u-s_2h)^t+\frac
1{n+1}\biggr)
 -C_{11}(u^{1-t}, u^t)\biggr|\to0,\cr\displaystyle
\sup_{a_n\le u\le1-b_n, -1\le s_1,s_2\le1}\biggl|C_{12}\biggl(\frac
{n}{n+1}(u-s_1h)^{1-t}+\frac1{n+1},\cr\displaystyle \hphantom{\sup_{a_n\le u\le1-b_n, -1\le s_1,s_2\le1}\biggl|C_{22}\biggl(}\frac{n}{n+1}(u-s_2h)^t+\frac
1{n+1}\biggr)
 -C_{12}(u^{1-t}, u^t)\biggr|\to0,\cr\displaystyle
\sup_{a_n\le u\le1-b_n, -1\le s_1,s_2\le1}\biggl|C_{22}\biggl(\frac
{n}{n+1}(u-s_1h)^{1-t}+\frac1{n+1},
\cr\displaystyle \hphantom{\sup_{a_n\le u\le1-b_n, -1\le s_1,s_2\le1}\biggl|C_{22}\biggl(}
\frac{n}{n+1}(u-s_2h)^t
+\frac1{n+1}\biggr)
 -C_{22}(u^{1-t}, u^t)\biggr|\to0.
}
\end{equation}
Thus, by (\ref{Th2-1}), (\ref{pfL1-2}), and similar arguments used in
the proof of Theorem~\ref{asm}, we can show that
%
\begin{eqnarray}
\label{pfL1-3}
\hspace*{-5pt}\sqrt n \mathit{II}_1\stackrel{d}{\to}
\int_0^1\!\{W(u^{1-t},u^t)-W(u^{1-t},1)C_1(u^{1-t},u^t)
-W(1,u^t)C_2(u^{1-t},u^t)\}\lambda(u,t)\,\mathrm{d}u.\qquad\quad
\end{eqnarray}

It is straightforward to verify that
%
\begin{equation}
\label{pfL1-4}
\cases{\displaystyle
|C_1(u^{1-t},u^t)u^{1-t}|=\mathrm{O}\bigl(u^{A(t)}\bigr)=\mathrm{O}(u^{1/2}),\vspace*{2pt}\cr\displaystyle
|C_2(u^{1-t},u^t)u^t|=\mathrm{O}\bigl(u^{A(t)}\bigr)=\mathrm{O}(u^{1/2}),\vspace*{2pt}\cr\displaystyle
|C_{11}(u^{1-t},u^t)u^{2-2t}\{1-\log u\}|=\mathrm{O}\bigl(u^{A(t)}\bigr)=\mathrm{O}(u^{1/2}),\vspace*{2pt}\cr\displaystyle
|C_{22}(u^{1-t},u^t)u^{2t}\log u|=\mathrm{O}\bigl(u^{A(t)}\bigr)=\mathrm{O}(u^{1/2}),\vspace*{2pt}\cr\displaystyle
|C_{12}(u^{1-t},u^t)u\{1-\log u\}|=\mathrm{O}\bigl(u^{A(t)}\bigr)=\mathrm{O}(u^{1/2})
}
\end{equation}
uniformly for $u\in[a_n, 1-b_n]$. By Taylor's
expansion, we have
%
\begin{eqnarray}
\label{II-2}
\hspace*{-15pt}\mathit{II}_2&=&\int_{a_n}^{1-b_n}\!\int_{-1}^1\!\int_{-1}^1\!
\biggl\{\!C_1(u^{1-t},u^t)u^{1-t} \biggl(\frac{n}{n+1}\biggl(1-\frac{s_1h}u\biggr)^{1-t}\!+\!\frac
1{(n+1)u^{1-t}}-1 \biggr)\nonumber\\
&&\hphantom{\int_{a_n}^{1-b_n}\int_{-1}^1\int_{-1}^1
\biggl\{}{} +C_2(u^{1-t},u^t)u^t \biggl(\frac{n}{n+1}\biggl(1-\frac{s_2h}u\biggr)^t+\frac
1{(n+1)u^{t}}-1 \biggr)\nonumber\\
&&\hphantom{\int_{a_n}^{1-b_n}\int_{-1}^1\int_{-1}^1
\biggl\{} {}+\frac12
C_{11}(u^{1-t},u^t)\bigl(1+\mathrm{o}(1)\bigr)u^{2-2t}\nonumber \\
&&\hphantom{\int_{a_n}^{1-b_n}\int_{-1}^1\int_{-1}^1
\biggl\{{}+\;}
{} \times\biggl(\frac{n}{n+1}\biggl(1-\frac
{s_1h}u\biggr)^{1-t}+\frac1{(n+1)u^{1-t}}-1 \biggr)^2\\
&&\hphantom{\int_{a_n}^{1-b_n}\int_{-1}^1\int_{-1}^1
\biggl\{}{} +\frac12C_{22}(u^{1-t},u^t)\bigl(1+\mathrm{o}(1)\bigr)u^{2t}\biggl (\frac n{n+1} \biggl(1-\frac
{s_2h}u\biggr)^t+\frac{1}{(n+1)u^t}-1 \biggr)^2
\nonumber\\
&&\hphantom{\int_{a_n}^{1-b_n}\int_{-1}^1\int_{-1}^1
\biggl\{}{} +C_{12}(u^{1-t},u^t)\bigl(1+\mathrm{o}(1)\bigr)
u \biggl(\frac n{n+1}\biggl(1-\frac{s_1h}u\biggr)^{1-t}+\frac
1{(n+1)u^{1-t}}-1 \biggr)\nonumber\\
&&\hphantom{\int_{a_n}^{1-b_n}\int_{-1}^1\int_{-1}^1
\biggl\{{}+\;}
{} \times\biggl(\frac{n}{n+1}\biggl(1-\frac{s_2h}u\biggr)^t+\frac
1{(n+1)u^t}-1 \biggr) \biggr\}
\nonumber\\
&&\hphantom{\int_{a_n}^{1-b_n}\int_{-1}^1\int_{-1}^1}
{} \times k(s_1)k(s_2)\lambda(u,t)\,\mathrm{d}s_1\,\mathrm{d}s_2\,\mathrm{d}u.
\nonumber
\end{eqnarray}
Consider the first term in the foregoing expression. By
\eqref{Th2-1}, \eqref{pfL1-2}, \eqref{pfL1-4}, and the symmetry of
$k(s)$, we have
\begin{eqnarray*}
&&\int_{a_n}^{1-b_n}\int_{-1}^1\int_{-1}^1
C_1(u^{1-t},u^t)u^{1-t} \biggl(\frac{n}{n+1}\biggl(1-\frac{s_1h}u\biggr)^{1-t}+\frac
1{(n+1)u^{1-t}}-1 \biggr)\\
&&\hphantom{\int_{a_n}^{1-b_n}\int_{-1}^1\int_{-1}^1}{}\times k(s_1)k(s_2)\lambda(u,t)\,\mathrm{d}s_1\,\mathrm{d}s_2\,\mathrm{d}u\\
&& \quad =\int_{a_n}^{1-b_n}\int_{-1}^1
C_1(u^{1-t},u^t)u^{1-t} \biggl(\frac{n}{n+1}\biggl(1-\frac{s_1h}u\biggr)^{1-t}-1
\biggr)k(s_1)\lambda(u,t)\,\mathrm{d}s_1\,\mathrm{d}u\\
&& \qquad {}+\frac{1}{n+1}\int_{a_n}^{1-b_n}C_1(u^{1-t},u^t)\lambda(u,t)\,\mathrm{d}u\\
&& \quad =\bigl(1+\mathrm{o}(1)\bigr)\int_{a_n}^{1-b_n}\int_{-1}^1
C_1(u^{1-t},u^t)u^{-1-t}\frac{nh^2}{2(n+1)}(1-t)(-t)s_1^2k(s_1)\lambda
(u,t)\,\mathrm{d}s_1\,\mathrm{d}u\\
&& \qquad {}+\frac{1}{n+1}\int_{a_n}^{1-b_n}C_1(u^{1-t},u^t)(1-u^{1-t})\lambda
(u,t)\,\mathrm{d}u\\
&& \quad =\mathrm{O}\biggl(h^2\int^{1-b_n}_{a_n}u^{-3/2}\lambda(u,t)\,\mathrm{d}u\biggr)+\mathrm{O}\biggl(n^{-1}\int
^{1-b_n}_{a_n}u^{-1/2}\lambda(u,t)\,\mathrm{d}u\biggr)=\mathrm{o}\bigl(1/\sqrt
n\bigr).
\end{eqnarray*}
Other terms of \eqref{II-2} can be handled in the same way, resulting in
%
\begin{eqnarray}
\label{pfL1-5}
\mathit{II}_2&=&\mathrm{o}\bigl(1/\sqrt
n\bigr)+\mathrm{O}\biggl(\int_{a_n}^{1-b_n}\biggl|C_2(u^{1-t},u^t)u^t\biggl(\frac{h^2}{u^2}+\frac
1{(n+1)u}\biggr)\lambda(u,t)\biggr|\,\mathrm{d}u \biggr)\qquad\nonumber\\
&&{} +\mathrm{O}\biggl(\int_{a_n}^{1-b_n}\biggl|C_{11}(u^{1-t},u^t)u^{2-2t}\biggl(\frac{h}{u}+\frac
1{(n+1)u}\biggr)^2\lambda(u,t)\biggr|\,\mathrm{d}u \biggr)\qquad\nonumber\\
&&{} +\mathrm{O}\biggl(\int_{a_n}^{1-b_n}\biggl|C_{22}(u^{1-t},u^t)u^{2t}\biggl(\frac{h}{u}+\frac
1{(n+1)u}\biggr)^2\lambda(u,t)\biggr|\,\mathrm{d}u \biggr)\qquad\nonumber\\
&&{} +\mathrm{O}\biggl(\int_{a_n}^{1-b_n}\biggl|C_{12}(u^{1-t},u^t)u\biggl(\frac{h}u+\frac
1{(n+1)u}\biggr)\biggl(\frac{h}u+\frac1{(n+1)u}\biggr)\lambda(u,t)\biggr|\,\mathrm{d}u \biggr) \qquad\ \ \\
&=&\mathrm{o}\bigl(1/\sqrt
n\bigr)+\mathrm{O}\biggl(h^2\int_{a_n}^{1-b_n}u^{-3/2}\lambda(u,t)\,\mathrm{d}u \biggr)\qquad\nonumber\\
&&{} +\mathrm{O}\biggl(h^2\int
_{a_n}^{1-b_n}\{\log u\}^{-1}u^{-3/2}\lambda(u,t)\,\mathrm{d}u \biggr)\qquad\nonumber\\
&=&\mathrm{o}\bigl(1/\sqrt n\bigr).\qquad\nonumber
\end{eqnarray}

For the second term, $I_2$, in (\ref{I1-4}), by the mean value theorem,
we can write
%
\begin{eqnarray}
\label{pfL1-a}
\hspace*{-10pt} I_2   &=&n^{-1}\int_{a_n}^{1-b_n}\sum_{i=1}^n\sum_{j=1}^n
\biggl\{\frac{\hat F_{n1,-i}^{1/(1-t)}(X_{j1})-\hat
F_{n1}^{1/(1-t)}(X_{j1})}{h} k\biggl(\frac{u-\hat
F_{n1}^{1/(1-t)}(X_{j1})}h\biggr)  \quad \qquad \nonumber\\
&& \hphantom{n^{-1}\int_{a_n}^{1-b_n}\sum_{i=1}^n\sum_{j=1}^n
\biggl\{}  {}+\frac12\biggl(\frac{\hat F_{n1,-i}^{1/(1-t)}(X_{j1})-\hat
F_{n1}^{1/(1-t)}(X_{j1})}{h}\biggr)^2
k'\biggl(\frac{u-\xi_{n,i,j}^{1/(1-t)}}h\biggr) \biggr\}  \quad \qquad  \\
&&  \hphantom{n^{-1}\int_{a_n}^{1-b_n}\sum_{i=1}^n\sum_{j=1}^n}{}\times K\biggl(\frac
{u-F_{n2}^{1/t}(X_{j2})}h\biggr)\lambda(u,t)\,\mathrm{d}u,  \quad \qquad
\nonumber
\end{eqnarray}
where $\xi_{n,i,j}$ is between $\hat F_{n1}(X_{j1})$ and $\hat
F_{n1,-i}(X_{j1})$. Using the equation
\[
F_{n1,-i}(X_{j1})-\hat F_{n1}(X_{j1})=\frac{1}{n}\hat
F_{n1}(X_{j1})-\frac{1}{n}I(X_{i1}\leq
X_{j1}),
\]
we have
%
\begin{equation}
\label{pfL1-6}
\cases{\displaystyle
\sup_{1\le i,j\le n}|\hat F_{n1}(X_{j1})-\hat F_{n1,-i}(X_{j1})|\le
n^{-1},\cr\displaystyle
\sup_{1\le i,j\le n}\bigl|\hat F_{n1}^{1/(1-t)}(X_{j1})-\hat
F_{n1,-i}^{1/(1-t)}(X_{j1})\bigr|\le
\frac{1}{1-t}n^{-1}.
}
\end{equation}
Then, uniformly for $u\in[a_n,1-b_n]$,
%
\begin{eqnarray}
\label{pfL1-6(a)}
&&\sum_{i=1}^n\sum_{j=1}^nI\biggl(\biggl|\frac{u-\xi_{n,i,j}^{1/(1-t)}}h\biggr|\le1\biggr) \qquad \nonumber\\
&& \quad \le \sum_{i=1}^n\sum_{j=1}^nP\biggl((u-h)^{1-t}-\frac1n\le\hat
F_{n1}(X_{j1})\le(u+h)^{1-t}+\frac1n\biggr) \qquad \nonumber
\\[-8pt]
\\[-8pt]
&& \quad \le n\times\biggl\{\frac{(n+1)(u+h)^{1-t}+(n+1)/n}{n}-\frac
{(n+1)(u-h)^{1-t}-(n+1)/n-1}n\biggr\} \qquad \nonumber\\
&& \quad = \mathrm{O}(u^{-1}nh) \qquad
\nonumber
\end{eqnarray}
and
%
\begin{equation}
\label{pfL1-6(b)}
\sum_{j=1}^nI\biggl(\biggl|\frac{u-F_{n1}^{1/(1-t)}(X_{j1})}h\biggr|\le1\biggr)=\mathrm{O}(u^{-1}h).
\end{equation}
Because $k(s)$ is a density function with support on $[-1, 1]$, it
follows from \eqref{pfL1-a}, \eqref{pfL1-6(a)} and
\eqref{pfL1-6(b)} that
%
\begin{eqnarray}
\label{pfL1-7}
I_2&=&\mathrm{O}\Biggl(h^{-1}n^{-2}\int_{a_n}^{1-b_n}
\sum_{i=1}^n\sum_{j=1}^nI\biggl(\biggl|\frac{u-F_{n1}^{1/(1-t)}(X_{j1})}h\biggr|\le
1\biggr)\lambda(u,t)\,\mathrm{d}u \Biggr)\nonumber\\
&&{} +\mathrm{O}\biggl(h^{-2}n^{-3}\int_{a_n}^{1-b_n}
\sum_{i=1}^n\sum_{j=1}^nI\biggl(\biggl|\frac{u-\xi_{n,i,j}^{1/(1-t)}}h\biggr|\le
1\biggr)\lambda(u,t)\,\mathrm{d}u \biggr)\\
&=&\mathrm{O}\biggl(n^{-1}\int_{a_n}^{1-b_n}u^{-1}\lambda(u,t)\,\mathrm{d}u \biggr)
=\mathrm{o}\bigl(1/\sqrt n\bigr).\nonumber
\end{eqnarray}
Similarly, we can show that
%
\begin{equation}
\label{pfL1-8}
I_3=\mathrm{o}\bigl(1/\sqrt n\bigr)\quad\mbox{and} \quad I_4=\mathrm{o}\bigl(1/\sqrt n\bigr).
\end{equation}
Thus, the lemma follows from (\ref{pfL1-3}), (\ref{pfL1-5}),
(\ref{pfL1-7}) and (\ref{pfL1-8}).
\end{pf}

%
\begin{lemma}\label{lem5.2} Under conditions of Theorem~\ref{jelm}, we have
\begin{eqnarray*}
\frac1n\sum_{i=1}^nQ_i^2(\theta_0)&\stackrel{p}{\to}& E \biggl(\int_0^1\{
W(u^{1-t},u^t)-W(u^{1-t},1)C_1(u^{1-t},u^t) \\
&& \hphantom{E \biggl(\int_0^1\{}{}-W(1,u^t)C_2(u^{1-t},u^t)\}\lambda(u,t)\,\mathrm{d}u \biggr)^2
\end{eqnarray*}
as $n\to\infty$.
\end{lemma}

\begin{pf} By \eqref{EXP}, we can write
\begin{eqnarray*}
&&Q_i^2(\theta)\\
&& \quad =\int_{a_n}^{1-b_n}\int_{a_n}^{1-b_n}\{\hat V_{i1}(u_1,t)\hat
V_{i1}(u_2,t)+\hat V_{i1}(u_1,t)\hat V_{i2}(u_2,t)-\hat
V_{i1}(u_1,t)u_2^{\theta}\\
&& \quad\hphantom{=\int_{a_n}^{1-b_n}\int_{a_n}^{1-b_n}\{} {}+\hat V_{i2}(u_1,t)\hat V_{i1}(u_2,t)+\hat V_{i2}(u_1,t)\hat
V_{i2}(u_2,t)-\hat V_{i2}(u_1,t)u_2^{\theta}\\
&& \quad \hphantom{=\int_{a_n}^{1-b_n}\int_{a_n}^{1-b_n}\{}{}
-u_1^{\theta}\hat V_{i1}(u_2,t)-u_1^{\theta}\hat
V_{i2}(u_2,t)+u_1^{\theta}u_2^{\theta}\}\lambda(u_1,t)\lambda(u_2,t)\,\mathrm{d}u_1\,\mathrm{d}u_2.
\end{eqnarray*}
Using arguments similar to those in (\ref{pfL1-a}), we have
%
\begin{eqnarray}
\label{pfL2-2a}
&&\frac1n\sum_{i=1}^n\int_{a_n}^{1-b_n}\int_{a_n}^{1-b_n}\hat
V_{i2}(u_1,t)\hat V_{i2}(u_2,t)\lambda(u_1,t)\lambda(u_2,t)\,\mathrm{d}u_1\,\mathrm{d}u_2\nonumber\\
&& \quad =\int_{a_n}^{1-b_n}\int_{a_n}^{1-b_n} \Biggl(\frac1n\sum_{i=1}^n\sum
_{j=1}^n\sum_{l=1}^n \biggl \{\frac{\hat F_{n1,-i}(X_{j1})-\hat F_{n1}(X_{j1})}h
\frac1{1-t}\hat F_{n1}^{t/(1-t)}(X_{j1})\nonumber\\
&& \quad \hphantom{=\int_{a_n}^{1-b_n}\int_{a_n}^{1-b_n} \Biggl(\frac1n\sum_{i=1}^n\sum
_{j=1}^n\sum_{l=1}^n \biggl \{}{}\times
k\biggl(\frac{u_1-\hat F_{n1}^{1/(1-t)}(X_{j1})}h\biggr)K\biggl(\frac{u_1-\hat
F_{n2}^{1/t}(X_{j2})}h\biggr)\nonumber\\
&& \quad \hphantom{=\int_{a_n}^{1-b_n}\int_{a_n}^{1-b_n} \Biggl(\frac1n\sum_{i=1}^n\sum
_{j=1}^n\sum_{l=1}^n \biggl \{}{}+\frac{\hat F_{n2,-i}(X_{j2})-\hat F_{n2}(X_{j2})}h\frac1t\hat
F_{n2}^{(1-t)/t}(X_{j2})\nonumber\\
&& \quad \hphantom{=\int_{a_n}^{1-b_n}\int_{a_n}^{1-b_n} \Biggl(\frac1n\sum_{i=1}^n\sum
_{j=1}^n\sum_{l=1}^n \biggl \{{}+\;}{}\times
k\biggl(\frac{u_1-\hat F_{n2}^{1/t}(X_{j2})}h\biggr)K\biggl(\frac{u_1-\hat
F_{n1}^{1/(1-t)}(X_{j1})}h\biggr) \biggr\}\nonumber
\\[-4pt]
\\[-8pt]
&& \quad \hphantom{=\int_{a_n}^{1-b_n}\int_{a_n}^{1-b_n} \Biggl(\frac1n\sum_{i=1}^n\sum
_{j=1}^n\sum_{l=1}^n}{} \times\biggl\{\frac{\hat F_{n1,-i}(X_{l1})-\hat F_{n1}(X_{l1})}h\frac1{1-t}\hat
F_{n1}^{t/(1-t)}(X_{l1})\nonumber\\
&& \quad \hphantom{=\int_{a_n}^{1-b_n}\int_{a_n}^{1-b_n} \Biggl(\frac1n\sum_{i=1}^n\sum
_{j=1}^n\sum_{l=1}^n \biggl \{{}\times\;}{}\times
k\biggl(\frac{u_2-\hat F_{n1}^{1/(1-t)}(X_{l1})}h\biggr)K\biggl(\frac{u_2-\hat
F_{n2}^{1/t}(X_{l2})}h\biggr)\nonumber\\
&& \quad \hphantom{=\int_{a_n}^{1-b_n}\int_{a_n}^{1-b_n} \Biggl(\frac1n\sum_{i=1}^n\sum
_{j=1}^n\sum_{l=1}^n \biggl \{{}\times\;}{} +\frac{\hat F_{n2,-i}(X_{l2})-\hat F_{n2}(X_{l2})}h
\frac1t\hat F_{n2}^{(1-t)/t}(X_{l2})\nonumber\\
\hspace*{-15pt}&& \quad \hphantom{=\int_{a_n}^{1-b_n}\int_{a_n}^{1-b_n} \Biggl(\frac1n\sum_{i=1}^n\sum
_{j=1}^n\sum_{l=1}^n \biggl \{{}\times\;{}+\;}{}\times
k\biggl(\frac{u_2-\hat F_{n2}^{1/t}(X_{l2})}h\biggr)
K\biggl(\frac{u_2-\hat F_{n1}^{1/(1-t)}(X_{l1})}h\biggr) \biggr\} \Biggr)\nonumber\quad\\
&& \quad \hphantom{=\int_{a_n}^{1-b_n}\int_{a_n}^{1-b_n}}{}\times\lambda(u_1,t)\lambda(u_2,t)\,\mathrm{d}u_1\,\mathrm{d}u_2+\mathrm{o}_p(1).\nonumber
\end{eqnarray}
It is straightforward to check that
\begin{eqnarray*}
&&\frac1n\sum_{i=1}^n\{\hat F_{n1}(x)-\hat F_{n1,-i}(x)\}\{\hat
F_{n1}(y)-\hat F_{n1,-i}(y)\}\\[-2pt]
&& \quad =\frac{n+1}{n^3}\hat F_{n1}(x\wedge y)-\frac{n+2}{n^3}\hat
F_{n1}(x)\hat F_{n1}(y),
\\[-2pt]
&&\frac1n\sum_{i=1}^n\{\hat F_{n2}(x)-\hat F_{n2,-i}(x)\}\{\hat
F_{n2}(y)-\hat F_{n2,-i}(y)\}\\[-2pt]
&& \quad =\frac{n+1}{n^3}\hat F_{n2}(x\wedge y)-\frac{n+2}{n^3}\hat
F_{n2}(x)\hat F_{n2}(y)
\end{eqnarray*}
and
\begin{eqnarray*}
&&\frac1n\sum_{i=1}^n\{\hat F_{n1}(x)-\hat F_{n1,-i}(x)\}\{\hat
F_{n2}(y)-\hat F_{n2,-i}(y)\}\\[-2pt]
&& \quad =\frac1{n^2}\hat C_n(\hat F_{n1}(x),\hat F_{n2}(y))-\frac
{n+2}{n^3}\hat F_{n1}(x)\hat F_{n2}(y).
\end{eqnarray*}
Then \eqref{pfL2-2a} can be written as
\begin{eqnarray*}
&&\frac1n\sum_{i=1}^n\int_{a_n}^{1-b_n}\int_{a_n}^{1-b_n}\hat
V_{i2}(u_1,t)\hat V_{i2}(u_2,t)\lambda(u_1,t)
\lambda(u_2,t)\,\mathrm{d}u_1\,\mathrm{d}u_2\\[-2pt]
&& \quad =\frac{1}{h^2}\int_{a_n}^{1-b_n}\int_{a_n}^{1-b_n} \Biggl(\frac
1{n^2h^2}\sum_{j=1}^n\sum_{l=1}^n   \biggl\{\bigl(\hat F_{n1}(X_{j1}\wedge X_{l1})-\hat F_{n1}(X_{j1})\hat
F_{n1}(X_{l1})\bigr)\frac1{(1-t)^2}\\[-2pt]
&& \quad \hphantom{=\frac{1}{h^2}\int_{a_n}^{1-b_n}\int_{a_n}^{1-b_n} \Biggl(\frac
1{n^2h^2}\sum_{j=1}^n\sum_{l=1}^n   \biggl\{}{}\times\hat F_{n1}^{t/(1-t)}(X_{j1})\hat F_{n1}^{t/(1-t)}(X_{l1})k\biggl(\frac
{u_1-\hat F_{n1}^{1/(1-t)}(X_{j1})}h\biggr)\\[-2pt]
&& \quad \hphantom{=\frac{1}{h^2}\int_{a_n}^{1-b_n}\int_{a_n}^{1-b_n} \Biggl(\frac
1{n^2h^2}\sum_{j=1}^n\sum_{l=1}^n   \biggl\{}{}\times
K\biggl(\frac{u_1-\hat F_{n2}^{1/t}(X_{j2})}h\biggr)k\biggl(\frac{u_2-\hat
F_{n1}^{1/(1-t)}(X_{l1})}h\biggr)
\\[-2pt]
&& \quad \hphantom{=\frac{1}{h^2}\int_{a_n}^{1-b_n}\int_{a_n}^{1-b_n} \Biggl(\frac
1{n^2h^2}\sum_{j=1}^n\sum_{l=1}^n   \biggl\{}{}\times
K\biggl(\frac{u_2-\hat F_{n2}^{1/t}(X_{l2})}h\biggr)\\[-2pt]
&& \quad \hphantom{=\frac{1}{h^2}\int_{a_n}^{1-b_n}\int_{a_n}^{1-b_n} \Biggl(\frac
1{n^2h^2}\sum_{j=1}^n\sum_{l=1}^n   \biggl\{}{}
+\bigl(\hat F_{n2}(X_{j2}\wedge X_{l2})-\hat F_{n2}(X_{j2})\hat
F_{n2}(X_{l2})\bigr)\frac1{t^2}\\[-2pt]
&& \quad \hphantom{=\frac{1}{h^2}\int_{a_n}^{1-b_n}\int_{a_n}^{1-b_n} \Biggl(\frac
1{n^2h^2}\sum_{j=1}^n\sum_{l=1}^n   \biggl\{{}+\;}{}\times
\hat F_{n2}^{(1-t)/t}(X_{j2})\hat F_{n2}^{(1-t)/t}(X_{l2})k\biggl(\frac
{u_1-\hat F_{n2}^{1/t}(X_{j2})}h\biggr)\\[-2pt]
&& \quad \hphantom{=\frac{1}{h^2}\int_{a_n}^{1-b_n}\int_{a_n}^{1-b_n} \Biggl(\frac
1{n^2h^2}\sum_{j=1}^n\sum_{l=1}^n   \biggl\{{}+\;}{}\times
K\biggl(\frac{u_1-\hat F_{n1}^{1/(1-t)}(X_{j1})}h\biggr)k\biggl(\frac{u_2-\hat
F_{n2}^{1/t}(X_{l2})}h\biggr)\\[-2pt]
&& \quad \hphantom{=\frac{1}{h^2}\int_{a_n}^{1-b_n}\int_{a_n}^{1-b_n} \Biggl(\frac
1{n^2h^2}\sum_{j=1}^n\sum_{l=1}^n   \biggl\{{}+\;}{}\times
K\biggl(\frac{u_2-\hat F_{n1}^{1/(1-t)}(X_{l1})}h\biggr)\\[-2pt]
&& \quad \hphantom{=\frac{1}{h^2}\int_{a_n}^{1-b_n}\int_{a_n}^{1-b_n} \Biggl(\frac
1{n^2h^2}\sum_{j=1}^n\sum_{l=1}^n   \biggl\{}{}
+\bigl(\hat C_n(\hat F_{n1}(X_{j1}),
\hat F_{n2}(X_{l2}))-\hat F_{n1}(X_{j1})\hat F_{n2}(X_{l2})\bigr)\\[-2pt]
&& \quad \hphantom{=\frac{1}{h^2}\int_{a_n}^{1-b_n}\int_{a_n}^{1-b_n} \Biggl(\frac
1{n^2h^2}\sum_{j=1}^n\sum_{l=1}^n   \biggl\{{}+\;}{}\times
\frac
1{t(1-t)}\hat F_{n1}^{1/(1-t)}(X_{j1})\hat F_{n2}^{(1-t)/t}(X_{l2})\\[-2pt]
&& \quad \hphantom{=\frac{1}{h^2}\int_{a_n}^{1-b_n}\int_{a_n}^{1-b_n} \Biggl(\frac
1{n^2h^2}\sum_{j=1}^n\sum_{l=1}^n   \biggl\{{}+\;}{}\times
k\biggl(\frac
{u_1-\hat F_{n1}^{1/(1-t)}(X_{j1})}h\biggr)K\biggl(\frac{u_1-\hat F_{n2}^{1/t}(X_{j2})}h\biggr)\\[-2pt]
&& \quad \hphantom{=\frac{1}{h^2}\int_{a_n}^{1-b_n}\int_{a_n}^{1-b_n} \Biggl(\frac
1{n^2h^2}\sum_{j=1}^n\sum_{l=1}^n   \biggl\{{}+\;}{}\times
k\biggl(\frac{u_2-\hat F_{n2}^{1/t}(X_{l2})}h\biggr)K\biggl(\frac{u_2-\hat
F_{n1}^{1/(1-t)}(X_{l1})}h\biggr)\\[-2pt]
&& \quad \hphantom{=\frac{1}{h^2}\int_{a_n}^{1-b_n}\int_{a_n}^{1-b_n} \Biggl(\frac
1{n^2h^2}\sum_{j=1}^n\sum_{l=1}^n   \biggl\{}{}
+\bigl(\hat C_n(\hat F_{n1}(X_{l1}),
\hat F_{n2}(X_{j2}))-\hat F_{n1}(X_{l1})\hat F_{n2}(X_{j2})\bigr)\\[-2pt]
&& \quad \hphantom{=\frac{1}{h^2}\int_{a_n}^{1-b_n}\int_{a_n}^{1-b_n} \Biggl(\frac
1{n^2h^2}\sum_{j=1}^n\sum_{l=1}^n   \biggl\{{}+\;}{}\times\frac
1{t(1-t)}
\hat F_{n1}^{1/(1-t)}(X_{l1})
\hat F_{n2}^{(1-t)/t}(X_{j2})\\[-2pt]
&& \quad \hphantom{=\frac{1}{h^2}\int_{a_n}^{1-b_n}\int_{a_n}^{1-b_n} \Biggl(\frac
1{n^2h^2}\sum_{j=1}^n\sum_{l=1}^n   \biggl\{{}+\;}{}\times
k\biggl(\frac{u_2-\hat
F_{n1}^{1/(1-t)}(X_{l1})}h\biggr)K\biggl(\frac{u_2-\hat F_{n2}^{1/t}(X_{l2})}h\biggr)
\\[-2pt]
&& \quad \hphantom{=\frac{1}{h^2}\int_{a_n}^{1-b_n}\int_{a_n}^{1-b_n} \Biggl(\frac
1{n^2h^2}\sum_{j=1}^n\sum_{l=1}^n   \biggl\{{}+\;}{}\times
k\biggl(\frac{u_1-\hat
F_{n2}^{1/t}(X_{j2})}h\biggr)
K\biggl(\frac{u_1-\hat F_{n1}^{1/(1-t)}(X_{j1})}h\biggr) \biggr\} \Biggr)\\[-2pt]
&& \quad \hphantom{=\frac{1}{h^2}\int_{a_n}^{1-b_n}\int_{a_n}^{1-b_n}
}{}\times\lambda(u_1,t)\lambda(u_2,t)\,\mathrm{d}u_1\,\mathrm{d}u_2+\mathrm{o}_p(1).
\end{eqnarray*}
Based on the foregoing decomposition, we can show that
%
\begin{eqnarray}
\label{pfL2-2}
&&\frac1n\sum_{i=1}^n\int_{a_n}^{1-b_n}\int_{a_n}^{1-b_n}\hat
V_{i2}(u_1,t)\hat V_{i2}(u_2,t)\lambda(u_1,t)\lambda(u_2,t)\,\mathrm{d}u_1\,\mathrm{d}u_2\nonumber\\
&& \quad =\int_{0}^1\int_0^1 \bigl(\{u_1^{1-t}\wedge u_2^{1-t}-u_1^{1-t}u_2^{1-t}\}
C_1(u_1^{1-t},u_1^t)C_1(u_2^{1-t},u_2^t) \nonumber\\
&& \quad \hphantom{=\int_{0}^1\int_0^1 \bigl(}{}+\{u_1^t\wedge u_2^t-u_1^tu_2^t\}
C_2(u_1^{1-t},u_1^t)C_2(u_2^{1-t},u_2^t)\nonumber
\\[-8pt]
\\[-8pt]
&& \quad \hphantom{=\int_{0}^1\int_0^1 \bigl(}{}+\{C(u_1^{1-t}, u_2^{t})-u_1^{1-t}u_2^{t}\}
C_1(u_1^{1-t},u_1^t)C_2(u_2^{1-t},u_2^t)\nonumber\\
&& \quad \hphantom{=\int_{0}^1\int_0^1 \bigl(}{} +\{C(u_2^{1-t},u_1^t)-u_2^{1-t}u_1^t\}
C_1(u_2^{1-t},u_2^t)C_2(u_1^{1-t},u_1^t) \bigr)\nonumber\\
&& \quad \hphantom{=\int_{0}^1\int_0^1}{}\times\lambda(u_1,t)\lambda(u_2,t)\,\mathrm{d}u_1\,\mathrm{d}u_2+\mathrm{o}_p(1).\nonumber
\end{eqnarray}
Similarly, we have
%
\begin{eqnarray}
\label{pfL2-4}
&&\frac1n\sum_{i=1}^n\int_{a_n}^{1-b_n}\int_{a_n}^{1-b_n}\hat
V_{i1}(u_1,t)\hat V_{i2}(u_2,t)\,\mathrm{d}u_1\,\mathrm{d}u_2\nonumber\\
&& \quad =\int_0^1\int_0^1\{
C(u_1^{1-t},u_1^t)u_2^{1-t}C_1(u_2^{1-t},u_2^t)-C(u_1^{1-t}\wedge
u_2^{1-t},u_1^t)C_1(u_2^{1-t},u_2^t)\nonumber
\\[-8pt]
\\[-8pt]
&& \quad \hphantom{=\int_0^1\int_0^1\{}{}+C(u_1^{1-t},u_1^t)u_2^tC_2(u_2^{1-t},u_2^t)-C(u_1^{1-t},u_1^t\wedge
u_2^t)C_2(u_2^{1-t},u_2^{t})\}\nonumber\qquad\\
&& \quad \hphantom{=\int_0^1\int_0^1}{}\times\lambda(u_1,t)\lambda(u_2,t)\,\mathrm{d}u_1\,\mathrm{d}u_2+\mathrm{o}_p(1),
\nonumber
\\
\label{pfL2-3}
&&\frac1n\sum_{i=1}^n\int_{a_n}^{1-b_n}\int_{a_n}^{1-b_n}\hat
V_{i1}(u_1,t)\hat
V_{i1}(u_2,t)\lambda(u_1,t)\lambda(u_2,t)\,\mathrm{d}u_1\,\mathrm{d}u_2\nonumber
\\[-8pt]
\\[-8pt]
&& \quad =\int_0^1\int_0^1C(u_1^{1-t}\wedge u_2^{1-t},u_1^t\wedge
u_2^t)\lambda(u_1,t)\lambda(u_2,t)\,\mathrm{d}u_1\,\mathrm{d}u_2+\mathrm{o}_p(1),
\nonumber
\\
\label{pfL2-5}
&&\frac1n\sum_{i=1}^n\int_{a_n}^{1-b_n}\int_{a_n}^{1-b_n}\hat
V_{i1}(u_1,t)u_2^{\theta_0}\lambda(u_1,t)\lambda(u_2,t)\,\mathrm{d}u_1\,\mathrm{d}u_2\nonumber
\\[-8pt]
\\[-8pt]
&& \quad =\int_0^1\int_0^1C(u_1^{1-t},u_1^t)C(u_2^{1-t},u_2^t)\lambda
(u_1,t)\lambda(u_2,t)\,\mathrm{d}u_1\,\mathrm{d}u_2+\mathrm{o}_p(1),
\nonumber
\end{eqnarray}
and
%
\begin{equation}
\label{pfL2-6}
\frac1n\sum_{i=1}^n\int_{a_n}^{1-b_n}\int_{a_n}^{1-b_n}\hat
V_{i2}(u_1,t)u_2^{\theta_0}\lambda(u_1,t)\lambda(u_2,t)\,\mathrm{d}u_1\,\mathrm{d}u_2=\mathrm{o}_p(1).
\end{equation}
Thus, the lemma follows from (\ref{pfL2-2})--(\ref{pfL2-6}) and the
fact that
\begin{eqnarray*}
&&E \biggl(\int_0^1\{W(u^{1-t},u^t)-W(u^{1-t},1)C_1(u^{1-t},u^t)   -W(1,u^t)C_2(u^{1-t},u^t)\}\lambda(u,t)\,\mathrm{d}u \biggr)^2\\
&& \quad =\int_0^1\int_0^1 \bigl\{C(u_1^{1-t}\wedge u_2^{1-t}, u_1^t\wedge
u_2^t)-C(u_1^{1-t},u_1^t)C(u_2^{1-t},u_2^t)\\
&& \quad\hphantom{=\int_0^1\int_0^1 \bigl\{} {}-\bigl(C(u_1^{1-t}\wedge
u_2^{1-t},u_1^t)-C(u_1^{1-t},u_1^t)u_2^{1-t}\bigr)C_1(u_2^{1-t},u_2^t)\\
&& \quad\hphantom{=\int_0^1\int_0^1 \bigl\{} {}-\bigl(C(u_1^{1-t},u_1^{t}\wedge
u_2^t)-C(u_1^{1-t},u_1^t)u_2^t\bigr)C_2(u_2^{1-t},u_2^t)\\
&& \quad\hphantom{=\int_0^1\int_0^1 \bigl\{} {}-\bigl(C(u_1^{1-t}\wedge
u_2^{1-t},u_2^{1-t})-u_1^{1-t}C(u_2^{1-t},u_2^t)\bigr)C_1(u_1^{1-t},u_1^t)\\
&& \quad\hphantom{=\int_0^1\int_0^1 \bigl\{} {}+(u_1^{1-t}\wedge
u_2^{1-t}-u_1^{1-t}u_2^{1-t})C_1(u_1^{1-t},u_1^t)C_1(u_2^{1-t},u_2^t)\\
&& \quad\hphantom{=\int_0^1\int_0^1 \bigl\{} {}+\bigl(C(u_1^{1-t},u_2^t)-u_1^{1-t}u_2^t\bigr)C_1(u_1^{1-t},u_1^t)C_2(u^{1-t}_2,u_2^t)\\
&& \quad\hphantom{=\int_0^1\int_0^1 \bigl\{} {}-\bigl(C(u_2^{1-t},u_1^t\wedge
u_2^t)-u_1^tC(u_2^{1-t},u_2^t)\bigr)C_2(u_1^{1-t},u_1^t)\\
&& \quad\hphantom{=\int_0^1\int_0^1 \bigl\{} {}+\bigl(C(u_2^{1-t},u_1^t)-u_2^{1-t}u_1^t\bigr)C_2(u_1^{1-t},u_1^t)C_1(u_2^{1-t},u_2^t)\\
&& \quad\hphantom{=\int_0^1\int_0^1 \bigl\{} {}+(u_1^t\wedge
u_2^t-u_1^tu_2^t)C_2(u_1^{1-t},u_1^t)C_2(u_2^{1-t},u_2^t) \bigr\}\lambda
(u_1,t)\lambda(u_2,t)\,\mathrm{d}u_1\,\mathrm{d}u_2.
\end{eqnarray*}
\upqed
\end{pf}
\begin{pf*}{Proof of Theorem~\ref{jelm}} Using similar expansions as in
the proof of Lemma~\ref{lem5.1}, we can show that $\max_{1\le i\le
n}|Q_i(\theta_0)|=\mathrm{o}_p(n^{1/2})$. Thus, using Lemmas~\ref{lem5.1}
and~\ref{lem5.2}
and standard arguments in expanding the empirical likelihood ratio
(see, e.g., Owen
\cite{r9}), we obtain that as $n\rightarrow
\infty$,
\[
l(\theta_0)=\Biggl\{\sum_{i=1}^nQ_i(\theta_0)\Biggr\}^2\bigg/\sum_{i=1}^nQ_i^2(\theta
_0)+\mathrm{o}_p(1)\stackrel{d}{\to}\chi^2(1).
\]
\upqed
\end{pf*}

\section*{Acknowledgements}
We thank the editor, an associate
editor and two reviewers for their helpful comments. Peng's research
was supported by NSA Grant H98230-10-1-0170 and NSF Grant
DMS-10-05336. Qian's research was partly supported by the National
Natural Science Foundation of China (Grant 10971068). Yang's research was
partly supported by the National Basic Research Program (973
Program) of China (Grant 2007CB814905) and the Key Program of National
Natural Science Foundation of China (Grant 11131002).\looseness=1

%

\printhistory


\begin{thebibliography}{20}

\bibitem{bd}
%
\begin{barticle}[mr]
\bauthor{\bsnm{B{\"u}cher},~\bfnm{Axel}\binits{A.}} \AND
\bauthor{\bsnm{Dette},~\bfnm{Holger}\binits{H.}}
(\byear{2010}).
\btitle{A note on bootstrap approximations for the empirical copula process}.
\bjournal{Statist. Probab. Lett.}
\bvolume{80}
\bpages{1925--1932}.
\bid{doi={10.1016/j.spl.2010.08.021}, issn={0167-7152}, mr={2734261}}
\bptok{imsref}%
\end{barticle}
%
\endbibitem

\bibitem{BDV}
%
\begin{bmisc}[auto:STB|2012/01/31|14:46:44]
\bauthor{\bsnm{B{\"u}cher},~\bfnm{A.}\binits{A.}},
\bauthor{\bsnm{Dette},~\bfnm{H.}\binits{H.}} \AND
\bauthor{\bsnm{Volgushev},~\bfnm{S.}\binits{S.}}
(\byear{2011}).
\bhowpublished{New estimators of the Pickands dependence function and a
test for extreme-value dependence. \textit{Ann. Statist.} \textbf{39} 1963--2006.}
\bptok{imsref}%
\end{bmisc}
%
\endbibitem

\bibitem{r1}
%
\begin{barticle}[mr]
\bauthor{\bsnm{Cap{\'e}ra{\`a}},~\bfnm{P.}\binits{P.}},
\bauthor{\bsnm{Foug{\`e}res},~\bfnm{A.~L.}\binits{A.L.}} \AND
\bauthor{\bsnm{Genest},~\bfnm{C.}\binits{C.}}
(\byear{1997}).
\btitle{A nonparametric estimation procedure for bivariate extreme value
copulas}.
\bjournal{Biometrika}
\bvolume{84}
\bpages{567--577}.
\bid{doi={10.1093/biomet/84.3.567}, issn={0006-3444}, mr={1603985}}
\bptok{imsref}%
\end{barticle}
%
\endbibitem

\bibitem{cpz}
%
\begin{barticle}[mr]
\bauthor{\bsnm{Chen},~\bfnm{Jian}\binits{J.}},
\bauthor{\bsnm{Peng},~\bfnm{Liang}\binits{L.}} \AND
\bauthor{\bsnm{Zhao},~\bfnm{Yichuan}\binits{Y.}}
(\byear{2009}).
\btitle{Empirical likelihood based confidence intervals for copulas}.
\bjournal{J. Multivariate Anal.}
\bvolume{100}
\bpages{137--151}.
\bid{doi={10.1016/j.jmva.2008.04.005}, issn={0047-259X}, mr={2460483}}
\bptok{imsref}%
\end{barticle}
%
\endbibitem

\bibitem{r2}
%
\begin{barticle}[mr]
\bauthor{\bsnm{Deheuvels},~\bfnm{Paul}\binits{P.}}
(\byear{1991}).
\btitle{On the limiting behavior of the {P}ickands estimator for bivariate
extreme-value distributions}.
\bjournal{Statist. Probab. Lett.}
\bvolume{12}
\bpages{429--439}.
\bid{doi={10.1016/0167-7152(91)90032-M}, issn={0167-7152}, mr={1142097}}
\bptok{imsref}%
\end{barticle}
%
\endbibitem

\bibitem{r3}
%
\begin{barticle}[mr]
\bauthor{\bsnm{Falk},~\bfnm{Michael}\binits{M.}} \AND
\bauthor{\bsnm{Reiss},~\bfnm{Rolf-Dieter}\binits{R.D.}}
(\byear{2005}).
\btitle{On {P}ickands coordinates in arbitrary dimensions}.
\bjournal{J.~Multivariate Anal.}
\bvolume{92}
\bpages{426--453}.
\bid{doi={10.1016/j.jmva.2003.10.006}, issn={0047-259X}, mr={2107885}}
\bptok{imsref}%
\end{barticle}
%
\endbibitem

\bibitem{r4}
%
\begin{barticle}[mr]
\bauthor{\bsnm{Fermanian},~\bfnm{Jean-David}\binits{J.D.}},
\bauthor{\bsnm{Radulovi{\'c}},~\bfnm{Dragan}\binits{D.}} \AND
\bauthor{\bsnm{Wegkamp},~\bfnm{Marten}\binits{M.}}
(\byear{2004}).
\btitle{Weak convergence of empirical copula processes}.
\bjournal{Bernoulli}
\bvolume{10}
\bpages{847--860}.
\bid{doi={10.3150/bj/1099579158}, issn={1350-7265}, mr={2093613}}
\bptok{imsref}%
\end{barticle}
%
\endbibitem

\bibitem{r5}
%
\begin{barticle}[mr]
\bauthor{\bsnm{Genest},~\bfnm{Christian}\binits{C.}} \AND
\bauthor{\bsnm{Segers},~\bfnm{Johan}\binits{J.}}
(\byear{2009}).
\btitle{Rank-based inference for bivariate extreme-value copulas}.
\bjournal{Ann. Statist.}
\bvolume{37}
\bpages{2990--3022}.
\bid{doi={10.1214/08-AOS672}, issn={0090-5364}, mr={2541453}}
\bptok{imsref}%
\end{barticle}
%
\endbibitem

\bibitem{r6}
%
\begin{barticle}[mr]
\bauthor{\bsnm{Gong},~\bfnm{Yun}\binits{Y.}},
\bauthor{\bsnm{Peng},~\bfnm{Liang}\binits{L.}} \AND
\bauthor{\bsnm{Qi},~\bfnm{Yongcheng}\binits{Y.}}
(\byear{2010}).
\btitle{Smoothed jackknife empirical likelihood method for {ROC} curve}.
\bjournal{J. Multivariate Anal.}
\bvolume{101}
\bpages{1520--1531}.
\bid{doi={10.1016/j.jmva.2010.01.012}, issn={0047-259X}, mr={2609511}}
\bptok{imsref}%
\end{barticle}
%
\endbibitem

\bibitem{r7}
%
\begin{barticle}[mr]
\bauthor{\bsnm{Hall},~\bfnm{Peter}\binits{P.}} \AND
\bauthor{\bsnm{Tajvidi},~\bfnm{Nader}\binits{N.}}
(\byear{2000}).
\btitle{Distribution and dependence-function estimation for bivariate
extreme-value distributions}.
\bjournal{Bernoulli}
\bvolume{6}
\bpages{835--844}.
\bid{doi={10.2307/3318758}, issn={1350-7265}, mr={1791904}}
\bptok{imsref}%
\end{barticle}
%
\endbibitem

\bibitem{r8}
%
\begin{barticle}[mr]
\bauthor{\bsnm{Jing},~\bfnm{Bing-Yi}\binits{B.Y.}},
\bauthor{\bsnm{Yuan},~\bfnm{Junqing}\binits{J.}} \AND
\bauthor{\bsnm{Zhou},~\bfnm{Wang}\binits{W.}}
(\byear{2009}).
\btitle{Jackknife empirical likelihood}.
\bjournal{J. Amer. Statist. Assoc.}
\bvolume{104}
\bpages{1224--1232}.
\bid{doi={10.1198/jasa.2009.tm08260}, issn={0162-1459}, mr={2562010}}
\bptok{imsref}%
\end{barticle}
%
\endbibitem

\bibitem{ky}
%
\begin{barticle}[mr]
\bauthor{\bsnm{Kojadinovic},~\bfnm{Ivan}\binits{I.}} \AND
\bauthor{\bsnm{Yan},~\bfnm{Jun}\binits{J.}}
(\byear{2010}).
\btitle{Nonparametric rank-based tests of bivariate extreme-value dependence}.
\bjournal{J. Multivariate Anal.}
\bvolume{101}
\bpages{2234--2249}.
\bid{doi={10.1016/j.jmva.2010.05.004}, issn={0047-259X}, mr={2671214}}
\bptok{imsref}%
\end{barticle}
%
\endbibitem

\bibitem{r9}
%
\begin{barticle}[mr]
\bauthor{\bsnm{Owen},~\bfnm{Art~B.}\binits{A.B.}}
(\byear{1988}).
\btitle{Empirical likelihood ratio confidence intervals for a single
functional}.
\bjournal{Biometrika}
\bvolume{75}
\bpages{237--249}.
\bid{doi={10.1093/biomet/75.2.237}, issn={0006-3444}, mr={0946049}}
\bptok{imsref}%
\end{barticle}
%
\endbibitem

\bibitem{r10}
%
\begin{barticle}[mr]
\bauthor{\bsnm{Owen},~\bfnm{Art}\binits{A.}}
(\byear{1990}).
\btitle{Empirical likelihood ratio confidence regions}.
\bjournal{Ann. Statist.}
\bvolume{18}
\bpages{90--120}.
\bid{doi={10.1214/aos/1176347494}, issn={0090-5364}, mr={1041387}}
\bptok{imsref}%
\end{barticle}
%
\endbibitem

\bibitem{r11}
%
\begin{bmisc}[auto:STB|2012/01/31|14:46:44]
\bauthor{\bsnm{Owen},~\bfnm{A.}\binits{A.}}
(\byear{2001}).
\bhowpublished{\textit{Empirical Likelihood.}
New York: Chapman \& Hall/CRC.}
\bptok{imsref}%
\end{bmisc}
%
\endbibitem

\bibitem{r12}
%
\begin{barticle}[mr]
\bauthor{\bsnm{Pickands},~\bfnm{James}\binits{J.} \bsuffix{III}}
(\byear{1981}).
\btitle{Multivariate extreme value distributions}.
\bjournal{Bull. Inst. Internat. Statist.}
\bvolume{49}
\bpages{859--878}.
\bid{mr={0820979}}
\bptnote{check related}
\bptok{imsref}%
\end{barticle}
%
\endbibitem

\bibitem{r13}
%
\begin{barticle}[mr]
\bauthor{\bsnm{Qin},~\bfnm{Jing}\binits{J.}} \AND
\bauthor{\bsnm{Lawless},~\bfnm{Jerry}\binits{J.}}
(\byear{1994}).
\btitle{Empirical likelihood and general estimating equations}.
\bjournal{Ann. Statist.}
\bvolume{22}
\bpages{300--325}.
\bid{doi={10.1214/aos/1176325370}, issn={0090-5364}, mr={1272085}}
\bptok{imsref}%
\end{barticle}
%
\endbibitem

\bibitem{seger}
%
\begin{bmisc}[auto:STB|2012/01/31|14:46:44]
\bauthor{\bsnm{Segers},~\bfnm{J.}\binits{J.}}
(\byear{2012}).
\bhowpublished{Asymptotics of empirical copula
processes under nonrestrictive smoothness assumptions.
\textit{Bernoulli} \textbf{18} 764--782.}
\bptok{imsref}%
\end{bmisc}
%
\endbibitem

\bibitem{z}
%
\begin{bmisc}[auto:STB|2012/01/31|14:46:44]
\bauthor{\bsnm{Zhou},~\bfnm{M.}\binits{M.}}
\bhowpublished{emplik: Empirical likelihood ratio for
censored/truncated data. R package version 0.9-3-1. Available at
\url{http://www.ms.uky.edu/\textasciitilde mai/splus/library/emplik/}}.
\bptok{imsref}%
\end{bmisc}
%
\endbibitem

\end{thebibliography}
\end{document}